\documentclass[12pt]{article}
\usepackage{geometry}                
\usepackage{graphicx}
\usepackage{amssymb}					
\usepackage{amsmath}					
\usepackage{mathrsfs}
\usepackage{tikz}
\usepackage{epstopdf}
\usepackage{amsthm}

\newtheorem*{theorem*}{Theorem}
\newtheorem*{lemma*}{Lemma}
\newtheorem*{example*}{Example}
\newtheorem{theorem}{Theorem}
\newtheorem{lemma}{Lemma}

\usepackage{atbegshi}
\AtBeginDocument{\AtBeginShipoutNext{\AtBeginShipoutDiscard}}
\addtocounter{page}{-1}

\def\b0{{\bf 0}}
\def\b1{{\bf 1}}

\def\n{\noindent}

\begin{document}
\title{On dispersability of some circulant graphs*} {\thanks{accepted for publication April 9, 2024 in {\it Journal of Graph Algorithms and Applications}}}

\author{
Paul C. Kainen\\
 \texttt{kainen@georgetown.edu}
\\
Samuel Joslin
\\
\texttt{ssj34@georgetown.edu}
\\
Shannon Overbay
\\
\texttt{overbay@gonzaga.edu}
}
\date{}                                           

\newcommand{\Addresses}{{
  \bigskip
  \footnotesize

\n
Paul C. Kainen, \textsc{Department of Mathematics and Statistics,\\ Georgetown University,
37th and O Streets, N.W., Washington DC 20057}\\
\\

\n

\par\nopagebreak
}}

\maketitle
\abstract{
\n 
The  matching book thickness of a graph is the least number of pages in a book embedding such that each page is a matching. A graph is  dispersable if its matching book thickness equals its maximum degree.  
Minimum page matching book embeddings are given for bipartite and for most non-bipartite circulants contained in the (Harary) cube of a cycle and for various higher-powers.
}

\smallskip

\n
{\bf Key Phrases}: {circulant graph, matching book thickness, dispersable or nearly dispersable graph, sparseness of a matching book embedding, polymerization}

\vspace{1 cm}
\n

\section{Introduction}

\label{sec:in}

Dispersable graphs were introduced in \cite{bk79}, where it was conjectured that all {\it bipartite regular graphs are dispersable}.  This was disproved by Alam et al. \cite{abgkp2018} who showed that the Gray and Folkman graphs, though regular bipartite, are not dispersable.  These counterexamples are edge-transitive but not vertex-transitive.  In \cite{abdgkp2021}, Alam et al. show existence of infinitely many counterexamples to \cite{bk79} and conjecture {\it bipartite vertex-transitive graphs are dispersable}.

In this paper, we consider families of circulant graphs. Circulants  \cite{harary-circ} are used in graph theory, computer science, network engineering, and dynamical systems (e.g., \cite{bermond,gomez, circ2022, hwang, strogatz}).  Terms are defined in the next section.

Matching book embeddings of bipartite circulants $C$ are given where the page number is equal to the vertex degree $\Delta(C)$, supporting the conjecture in \cite{abdgkp2021}.  It can be shown that regular dispersable graphs must be bipartite \cite{so-thesis}. A nonbipartite circulant is {\it nearly dispersable} if one extra page suffices \cite{so-thesis}.  So far,
all nonbipartite circulants have been nearly dispersable and we conjecture here that {\it  nonbipartite, vertex-transitive graphs are nearly dispersable}.

Previous results support both conjectures.  For the complete bipartite  graph and the hypercube, see \cite{bk79}; for complete graphs and other bipartite graphs, see \cite{so-thesis}.  Cartesian products of even cycles are dispersable; even times odd cycles are nearly dispersable \cite{pck-bica}; and short odd (length at most 5) and arbitrary odd cycles have nearly dispersable product \cite{jko-2021b}.  Other classes of vertex-transitive graph that are known to be nearly dispersable include  the product of two {\it arbitrary} cycles  and of cycles with complete graphs, see \cite{shao-et-al, shao-zeling}, and some products of bipartite and nonbipartite graphs \cite{pupyrev}.   See also \cite{yu-shao-li} and \S 8.  Some graphs which are {\it not} vertex transitive also are dispersable such as trees \cite{so-thesis}, Halin trees \cite{shao-3}, and cubic planar bipartite graphs \cite{abdgkp2021, ko-2020}.

To define good matching book embeddings for an infinite family of graphs, one needs to give both layout and coloring schemes: algorithms which produce the needed vertex order and edge-to-page assignment from the various integers that identify each graph in the family.  Most of our families consist of circulants $C(n,S)$ with a fixed jump-length set $S$ and with the number $n$ of vertices reduced modulo 2 or 4.  The coloring algorithms can either be static (as in tables based on modularity) or dynamic (as in prescriptions for Hamiltonian cycles or paths).  See proofs of Theorems \ref{th:n-1-3-even} and \ref{th:cn12}, resp.

It turns out, however, that for nonbipartite circulants, perfectly regular patterns almost never succeed and irregularity is forced. 
Irregular features appear in two different ways: local and global.  

The local type of exception is involved in the ``twist'' (see \S 3) while the global type manifests as ``sparseness'' in many  examples, where one nearly reaches the lower bound except for a ``sparse" page with a small and structurally defined set of exceptional edges.  Computer search \cite[p 7]{pck-circ} gives random vertex-order, while our vertex-orders and edge-to-page functions are quite regular. Nevertheless, the edges of the sparse page are irregularly distributed in a characteristic pattern for all parameter values with the same modularity.

Any strategy to achieve the minimum number of pages in a matching book embedding of a graph family, such as $C(n,S)$, based on regular layout and algorithms  is a kind of ``polymerization process'' since almost all edges are placed in a repeated pattern. 
A polymer is a molecule composed of a sequence of many parts such as proteins composed of amino acids or RNA/DNA as a sequence of nucleotides. The sequence of parts may form a path or a cycle.

Here are four examples of {\it generalized polymerization}, a phenomenon that we believe deserves more thorough investigation. In each of these examples, a finite set of adjustments permits regularity for the arbitrarily large remainder.

\n
(1) The coloring irregularities given in \cite{pck-circ} for $C(2k+r,\{1,k\})$, $r \in \{0,1,2\}$, with up to 5 edges on the sparse page, allow a minimum page embedding.

\n
(2) Layouts and page-partition of the Cartesian product of $C_3$ or $C_5$ with another cycle \cite{jko-2021b} use a ``seed'' that is an exceptional copy of a repeated motif.

\n
(3) The twist used in \S 3 below for the case $C(2k,\{1,3\})$ with a dispersable circulant and in \S 4 for $C(4k+3,\{1,3\})$ in the nearly dispersable case.

\n
(4) Most of the matching book embeddings in this paper have a sparse page.

In contrast, {\it strict polymerization} is defined below to be an algorithmic procedure which puts together certain modular units with no adjustments.

The paper is organized as follows: 
\S 2 has definitions,  \S3 shows $C(n,\{1,3\})$ is dispersable for $n$ even, while \S 4 shows $C(n, \{1,3\})$ is nearly dispersable (n.d.) for $n$ odd. In \S 5 and \S 6, we prove $C(n, \{1,2\})$ and $C(n,\{2,3\})$ are n.d., and \S 7 shows $C(n,\{1,2,3\})$ is n.d. when $n$ is odd or a multiple of $7$ or $12$ (so $\geq 64.3\,$\% of the $C(n,\{1,2,3\})$ circulants); \S 7 also shows that $K_{2k} - k K_2$ and $K_{2k+1} - C_{2k+1}$ are n.d.
The last section has applications and a  discussion. 

\section{Definitions}

Undefined terms are as in \cite{harary}.   

The {\bf circulant graph} $C(n,S)$ of order $n$ with {\bf jump} set 
$S = \{i_1,\ldots,i_k\}$ 
is the graph on $[n] := \{1,\ldots,n\}$, where $j \in [n]$ is adjacent to $j+i_r$ (addition mod $n$), $r = 1,\ldots, k$ and $1 \leq i_1 < i_2 < \cdots < i_k \leq \lfloor n/2 \rfloor$, $k \geq 1$.  A graph is {\bf vertex-transitive} if for any two vertices, there is an isomorphism carrying one to the other.  The $k$-th {\bf Harary power} $C_n^k$ of an order-$n$ cycle $C_n$ \cite[p 14]{harary} is the circulant $C(n,[k])$, and any circulant of order $n$ with maximum jump $k$ is a vertex-transitive subgraph of $C_n^k$. 
The {\bf cube} of a cycle is the 3rd power.

A drawing of a graph is {\bf outerplane} (or {\bf convex} or {\bf circular}) if its vertices are placed along a circle (or the boundary of any convex region) and the edges are straight lines.
Two edges in an outerplane drawing {\bf cross} if they intersect at a non-endpoint. Let $(G,\omega)$ denote the outerplane drawing of a graph $G$ with cyclic order $\omega$ on $V(G)$.

A {\bf book embedding} \cite{bk79} of a graph $G$ is an outerplane drawing and an edge-partition such that edges in the same part do not cross. The parts of the partition are the {\bf pages} of the book embedding.  The {\bf book thickness} $bt(G)$ of $G$ is the least number of pages in any book embedding while $bt(G,\omega)$ is the least number of pages for the outerplane drawing $(G,\omega)$.

A {\bf proper edge-coloring} $c$ of a graph $G$ is a function $c: E(G) \to \{1,\ldots,r\}$ (the set of colors) such that adjacent edges get different colors.  Let $\chi'(G)$ be the least number of colors in a proper edge-coloring.  The remarkable theorem of Vizing \cite[p 133]{harary} states that $\chi'(G) \in \{\Delta(G), 1+\Delta(G)\}$ for all graphs $G$.

A {\bf matching book embedding} is a book embedding where the pages are matchings (no two edges are adjacent).  The {\bf matching book thickness} of a graph is the least number of pages in any matching book embedding; we write $mbt(G)$ or $mbt(G,\omega)$ as for book thickness. If $c$ is the edge-coloring determined by the pages, then the matching book embedding is the triple $(G,\omega, c)$.

Clearly, for every graph $G$, we have $\Delta(G) \leq \chi'(G) \leq mbt(G)$. 
A matching book embedding $(G,\omega, c)$ is {\bf dispersable} if the number $|c|$ of colors equals $\Delta(G)$ and is {\bf nearly dispersable} \cite{jko-2021b} if $|c| = 1 +\Delta(G)$. A graph is dispersable if it has a dispersable embedding and is nearly dispersable if it is not dispersable and has a nearly dispersable embedding. 
If $G$ is regular and dispersable, then it is bipartite \cite{so-thesis}.  
The {\bf sparseness} $s(G,\omega,c)$ of a nearly dispersable book embedding is the least number of edges on any page. The sparseness $s(G)$ of a nearly dispersable graph $G$ is the minimum sparseness over all minimum-page matching book embeddings.  

\begin{lemma}
Let $G$ be a regular nearly dispersable graph of order $n$. Then the sparseness of $G$ is at least 1 if $n$ is even and at least $\Delta/2 $ if $n$ is odd.
\label{lm:sparse}
\end{lemma}
\begin{proof}
For $n$ even, this is in Overbay \cite{so-thesis}, while for $n$ odd, each page has at least one uncovered vertex, so for any set of $\Delta$ pages, there is a set of $\geq \Delta$ distinct points which need to be covered by edges from the remaining page.
\end{proof}
An infinite sequence $\{(G_n, \omega_n, c_n)\}_{n \geq 1}$ of matching book embeddings, such that $|E(G_n)|$ is strictly increasing, is called {\bf sparse} (with sparseness $s$) if there exists $k$ such that, for all $n,\;$ we have (i) $\Delta(G_n)=k$, (ii) $mbt(G_n, \omega_n)= k{+}1$, (iii) $c_n: E(G_n) \to [k{+}1]$ is onto, and (iv) $s(G_n, \omega_n, c_n)=s$.
For each $n$, the page with $s$ edges is the {\bf sparse} or ``exceptional''  page \cite{pck-circ};
cf. \cite{Duj2004}.

It remains for us to define the $m$-fold {\it polymerization} of a matching book embedding of a circulant to form a circulant with the same set of jump-lengths but $m$-fold more vertices with no increase in the number of pages.

For $n \geq 7$ and nonempty $S \subseteq \{1, 2, \ldots, \lceil \frac{n-1}{2} \rceil\}$, let $(C(n, S), \nu_n, c)$, $\nu_n := (1,\ldots,n)$, be a matching book embedding. We call an edge $e = a_i a_k \in E := E(C(n, S))$ {\bf long} if $d_C(a_i, a_k) < |k - i|$ and  {\bf short} if $d_C(a_i, a_k) = |k-i|$, where $d_C(u,w)$ denotes the $C_n$-distance between two vertices $u$ and $w$, where $C_n$ is the graph induced by the cyclic vertex order.   
The sets $E_\Lambda$ and $E_\Sigma$ of long and short edges form a nontrivial partition of $E$. If $1 \in S$, then the edges $a_i a_{i+1}$ are short for $i = 1, \ldots, n-1$ but the edge $a_1 a_n$ is long.  If $n=8$ and $3 \in S$, then $a_i a_{i+3} \in E_s$ for $1 \leq i \leq 5$, while $a_6 a_1$, $a_7 a_2$, and $a_8 a_3$ are long.

\begin{lemma}
Let $m \geq 2, n \geq 7$.  If $S \subseteq \{1, 2, \ldots, \lceil \frac{n-1}{2} \rceil\}$, then 
\begin{equation}
mbt(C(nm, S), \nu_{nm}) \leq mbt(C(n,S), \nu_n).
\label{eq:polym}
\end{equation}
\label{lm:polymer}
\end{lemma}
\begin{proof}
Place $m$ copies of the vertex set of $C(n,S)$ from left to right, where, for the $j$-th copy, the vertices $a_1^j, \ldots, a_n^j$ are placed from left to right.   Thus,
$$(a_1^1, \dots, a_n^1, a_1^2, \ldots, a_n^2, a_1^3, \ldots, a_n^m)$$
is a list of the $nm$ vertices in the $m$ copies of $C(n,S)$.
Put in all the short edges for all the copies. For $j = 1, \ldots, m-1$, each long edge $a_i^j a_k^j$ with $k > i$ is replaced by $a_k^j a_i^{j+1}$ and each long edge $a_i^m a_k^m$ is replaced by $a_k^m a_i^1$ if $k > i$.  One obtains $C(nm, S)$ with vertex order $\nu_{nm}$.  Use the same coloring $c$ for the edges in $C(nm, S)$ as in the matching book embedding for $C(n, S)$.  Long and short edges cross in the $j$-th copy if and only if their images under the edge-rearrangement cross correspondingly. Hence, $c$ is a page assignment.
\end{proof}
This process defines the $m$-fold {\bf strict polymerization} of the circulant and of its matching book embedding.  
Note that equality can fail to hold in (\ref{eq:polym}) - e.g., for an even polymerization of an odd cycle.
If $(C(n,S),\nu_n,c)$ is dispersable, then so is $C(nm, S), \nu_{nm}, c)$, while if $(C(n,S),\nu_n,c)$ is nearly dispersable, then $C(nm, S), \nu_{nm}, c)$ is either dispersable or nearly dispersable. If both 
are nearly dispersable, then
$s(C(nm, S), \nu_{nm}, c) \leq m \cdot s(C(n,S), \nu_n, c)$.

\section{The bipartite case of $C(n,\{1,3\})$}

In this section we show that the circulants $C(n, \{1,3\})$, $n \geq 6$ even, are dispersable.  When $n=6$, the corresponding circulant is $K_{3,3}$, and this graph, as with {\it all} complete regular bipartite graphs $ K_{a,a}$, is dispersable \cite{bk79}, \cite[p. 88]{so-thesis}.

\begin{figure}[ht!]
\centering
\includegraphics[width=50mm]{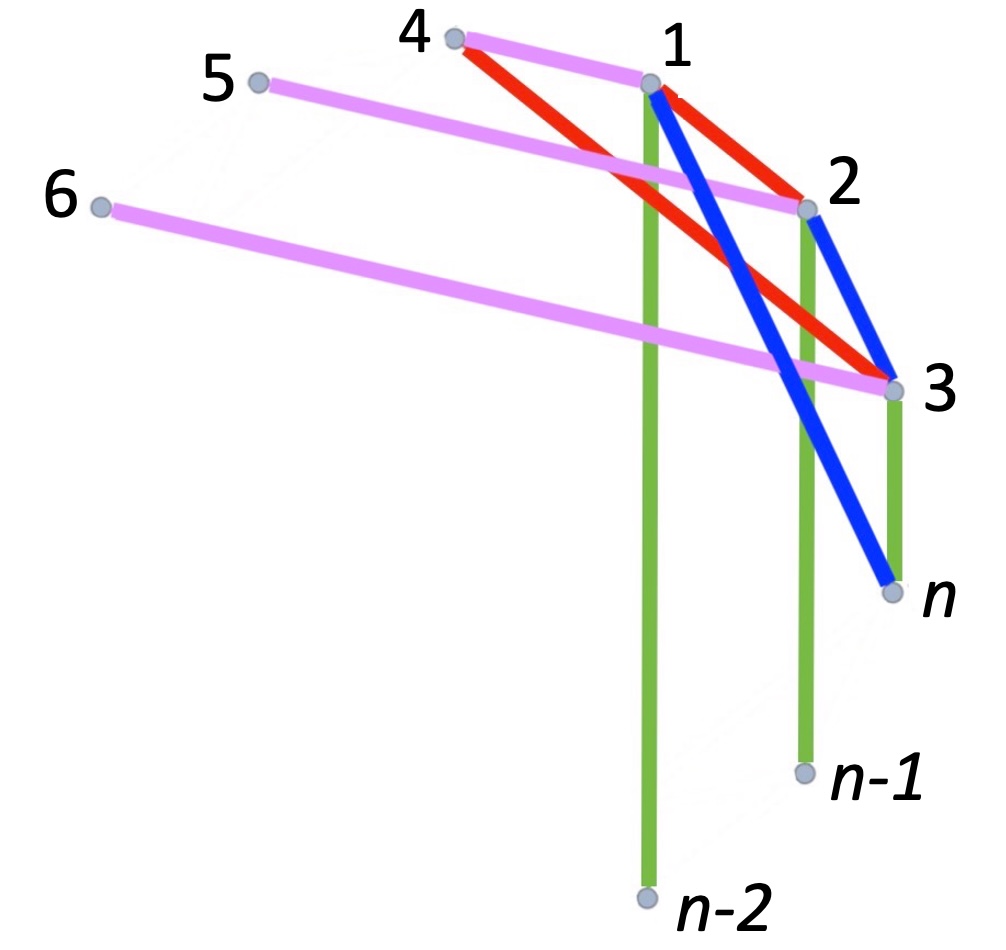}
\caption{Common four-coloring $c$ of the twist \label{Figure 1}}
\end{figure}

\begin{figure}[ht!]
\centering
\includegraphics[width=60mm]{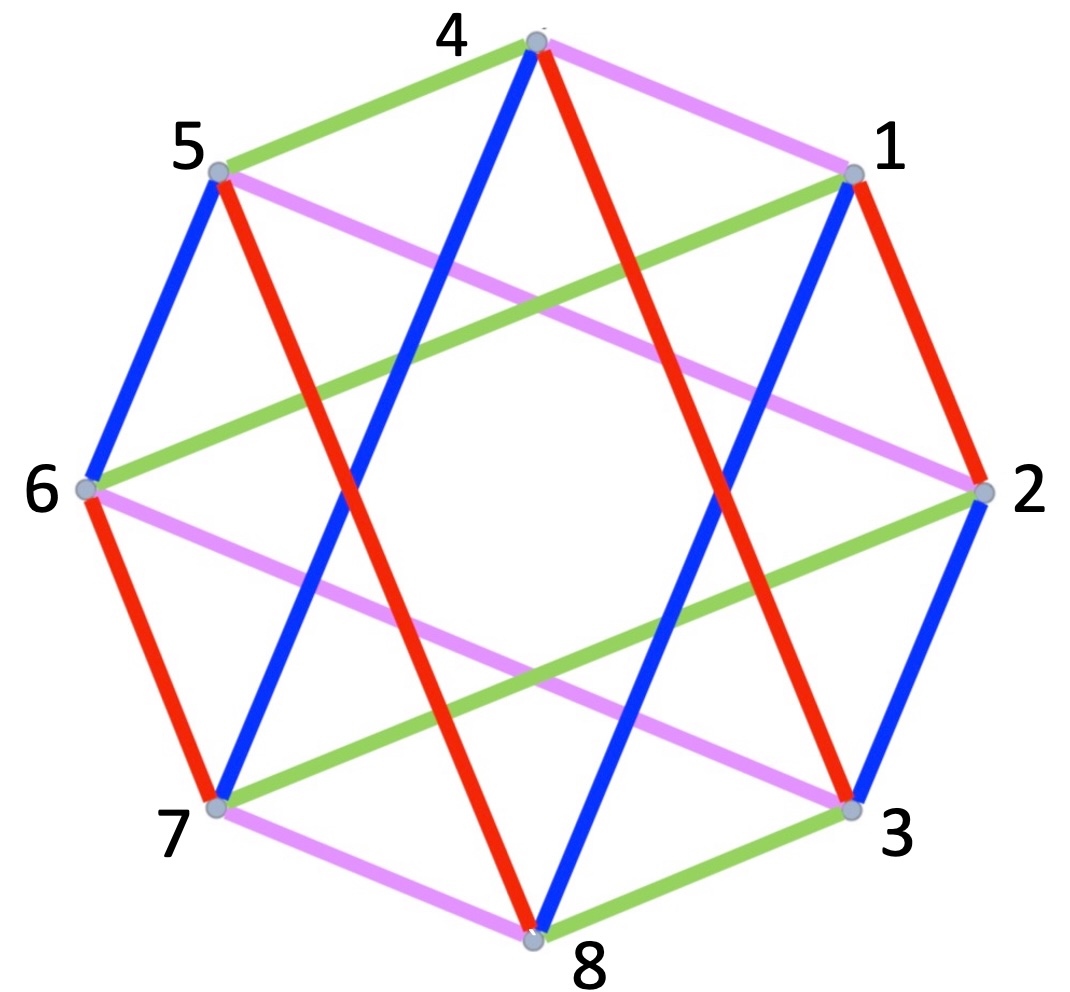}
\includegraphics[width=56mm]{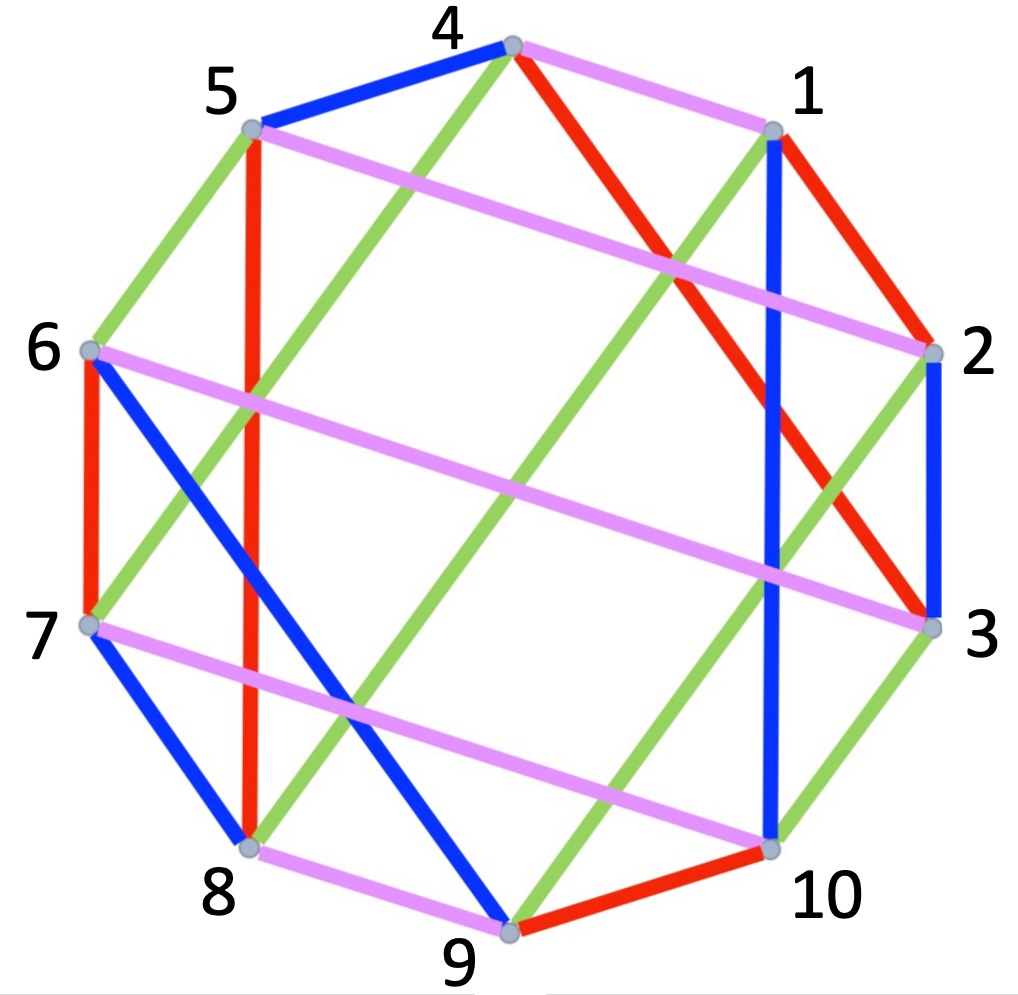}
\caption{Dispersable embeddings of $C(8,\{1,3\})$ and $C(10,\{1,3\})$ \label{Figure 2}}
\end{figure}

\begin{theorem} Let $n\geq 8$, $n$ even, and put $\tau_n := (1, 2, 3, n, n-1, n-2, \ldots, 6, 5, 4)$.
Then $(C(n, \{1,3\}), \tau_n, c_n)$
is a dispersable book embedding {\rm ($c_n$ defined below)}.
\label{th:n-1-3-even}
\end{theorem}

\begin{proof}
Let 
$c$ be the 4-coloring  of the 10 edges in the ``twist'' shown in Fig. 
\ref{Figure 1}: 
\begin{itemize}
\item Red: $1{-}2, \, 3{-}4$ (i.e., color of edge $1{{-}}2$ is red, etc.),
\item Blue: $2{-}3,\, 1{-}n$,
\item Green: $3{-}n,\, 2{-}(n{-}1), \, 1{-}(n{-}2)$,
\item Purple: $1{-}4,\, 2{-}5,\,3{-}6$.
\end{itemize}
 {\bf Case 1}:  
$n = 4k$, $k \geq 2$. Assign the non-twist edges to four pages as follows:
\begin{itemize}
\item Blue: $a{-}(a{+}3)$, $\;(a{+}1){-}(a{+}2)$,  $a = 4t$, $t \in [k-1]$, 
\item Red: $a{-}(a{+}3)$, $\;(a{+}1){-}(a{+}2)$, $a = 1+4t$, $t \in [k-1]$,
\item Green: $a{-}(a{+}3)$, $\;(a{+}1){-}(a{+}2)$, $a =2+4t$, $t \in [k-2]$, and $4{-}5$,
\item Purple:  $a{-}(a{+}3)$, $\;(a{+}1){-}(a{+}2)$,  $a = 3+4t$, $t \in [k-2]$, and $(n{-}1){-}(n)$.
\end{itemize}
All $8k$ edges of $C(n, \{1,3\})$ appear, $2k$ in a 
page. 
In Fig. 2, left, $k=2$, so $[k-2] = \emptyset$; hence, there is only one green edge not in the twist coloring.  

On the Red and Blue pages, one has 2 edges on the common twist and $k-1$ edges of types $a{-}(a+3)$ and another $k-1$ edges of type $(a{+}1){-}(a{+}2)$. 
Similarly, on the Purple and Green pages, one has 3 edges on the common twist, an additional edge ($4{-}5$ or $(n{-}1){-}n$) and $k-2$ each of types $a{-}(a{+}3)$ and $(a{+}1){-}(a{+}2)$. These $8k$ edges are distinct and exhaust the edges of $C(n, \{1,3\})$. 
For $n=8$, $k=2$ so $t=1$ and on the Red page, $a=5$. In Fig. \ref{Figure 2}, after
$1{-}2$ and $3{-}4$, we also have $5{-}8$ and $6{-}7$ on the red page.

By definition, the edges on each of these pages are pairwise-disjoint, while 
pages are crossing-free since the edges in a color class can only be (i) non-crossing edges of the common twist, (ii) isolated edges on the outer cycle of the form $4{-}5$ or $(n{-}1){-}n$, (iii) nested edges of the form $\{ a{-}(a{+}3), (a{+}1){-}(a{+}2) \}$.

\vskip 0.1in
\begin{figure}[ht!]
\centering
\includegraphics[width=60mm]{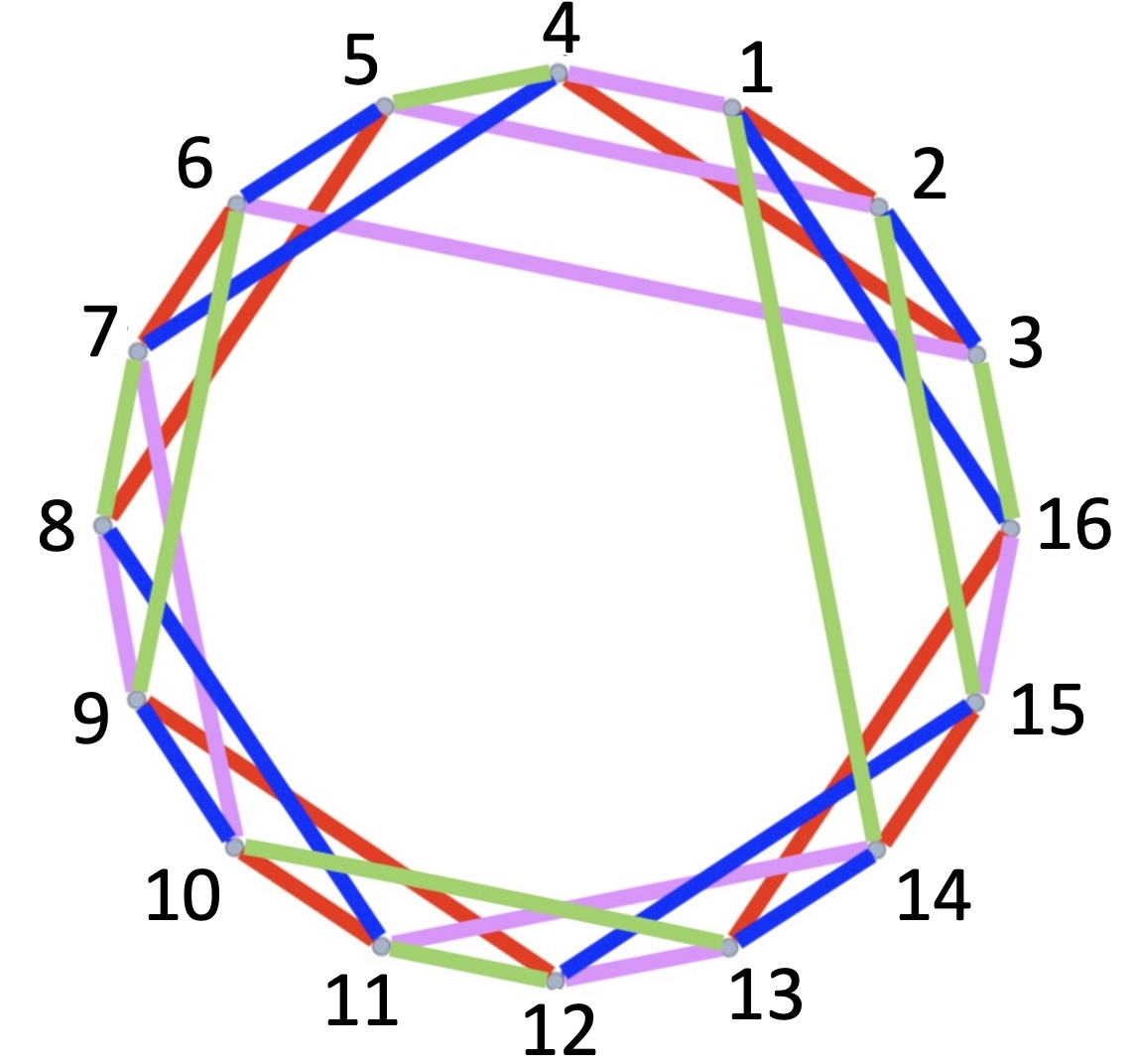}
\includegraphics[width=60mm]{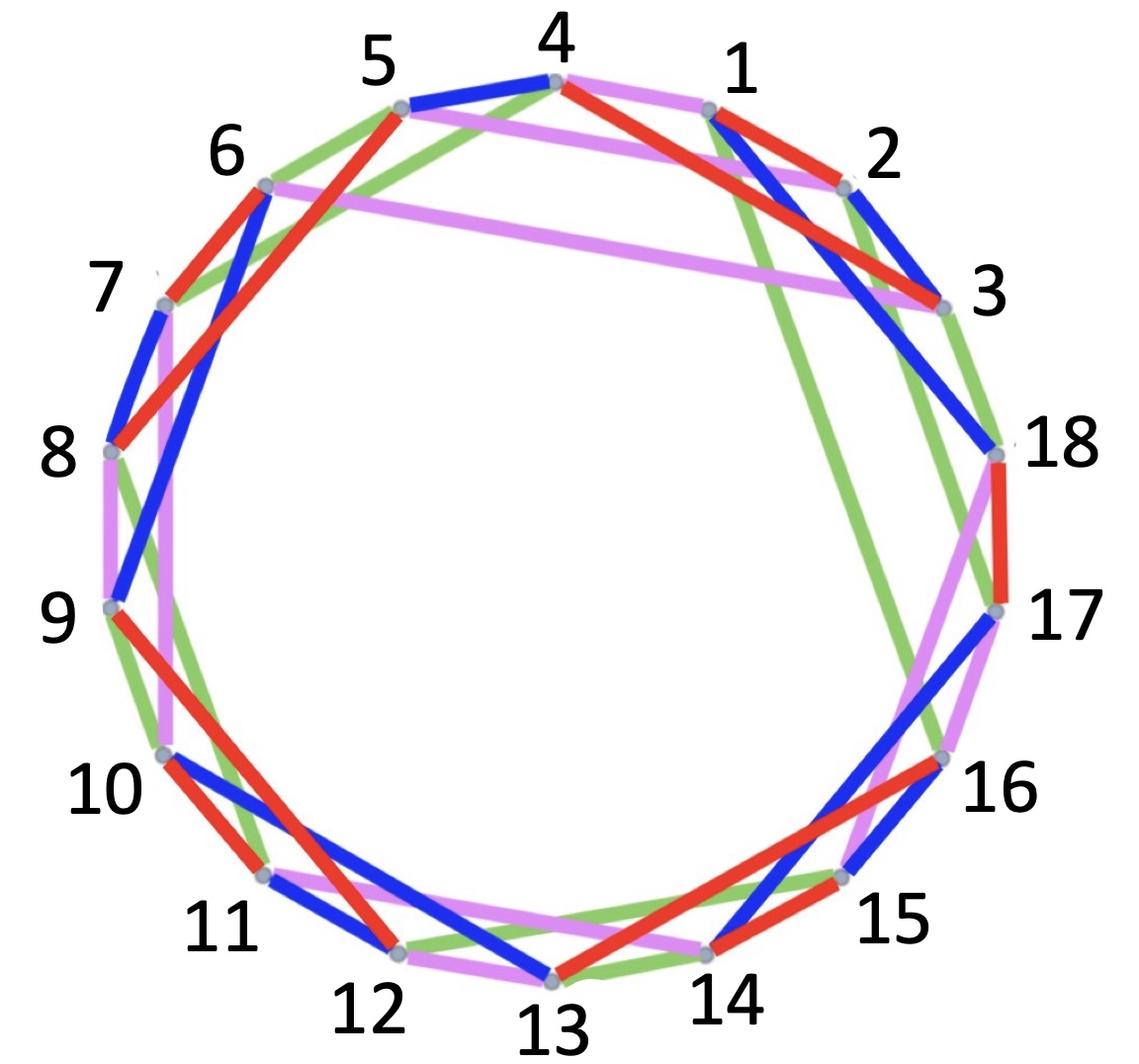}
\caption{Dispersable embeddings of $C(16,\{1,3\})$ and $C(18,\{1,3\})$\label{Figure 3}}
\end{figure}

{\bf Case 2}:	 Let $n = 4k + 2$,  $k \geq 2$; put non-twist edges in four pages as follows:
\begin{itemize}
\item Green: $a{-}(a{+}3)$, $(a{+}1){-}(a{+}2)$, $a = 4t$, $t \in [k-1]$,
\item Red: $a{-}(a{+}3)$, $(a{+}1){-}(a{+}2)$, $a = 1+4t$, $t \in [k-1]$, and $n{-}(n{-}1)$,
\item Blue: $a{-}(a{+}3)$, $(a{+}1){-}(a{+}2)$, $a = 2+4t$, $t \in [k-1]$, and $4{-}5$,
\item Purple:  $a{-}(a{+}3)$, $(a{+}1){-}(a{+}2)$, where $a = 3+4t$, $t \in [k-1]$.
\end{itemize}
See Fig. \ref{Figure 2}, right.  Fig. \ref{Figure 3} shows both schemes for $k=4$.
\end{proof}



\section{The nonbipartite case of $C(n, \{1,3\})$}

When $n$ is odd, $C(n, \{1,3\})$ is not bipartite. But it is nearly dispersable.

\begin{theorem}
Let $n \geq 7$ be odd.  Then $C(n,\{1,3\})$ is nearly dispersable. 
\label{th:cn13-odd}
\end{theorem}

\begin{proof}
{\bf Case 1}: $n = 4k+1$. Let $\nu_n$ be the natural ordering $(1, 2, 3, \ldots, n)$ around the circle. We show that $mbt(C(n,\{1,3\}), \nu_n) = 5$, which is minimum for a nonbipartite 4-regular graph.
The $2n = 8k + 2$ edges are of the form $u{-}(u{+}3)$ or $u{-}(u{+}1)$, $u \in [n]$, with addition modulo $n$.

Assign each of the four colors red, purple, green, and blue to $2k$ edges and the fifth color, black, to the two remaining edges as follows:
\begin{itemize}
\item Red: $\;\;\;\;a{-}(a{+}3)$, $(a{+}1){-}(a{+}2)$, $a = 1+4t$, $t \in [k]-1$,
\item Purple: $a{-}(a{+}3)$, $(a{+}1){-}(a{+}2)$, $a = 2+4t, \,t\in [k]-1$,
\item Green: $\;a{-}(a{+}3)$, $(a{+}1){-}(a{+}2)$, $a = 3+4t, \,t\in [k]-1$,
\item Blue: $\;\;\;a{-}(a{+}3)$, $(a{+}1){-}(a{+}2)$, $a = 4+4t$, $t \in [k]-1$,
\item Black: $1{-}2$, $3{-}n$ (sparseness is 2; as $\Delta = 4$, this is the minimum).
\end{itemize}

Since the edge pairs of the form $\{a{-}(a{+}3)$, $(a{+}1){-}(a{+}2)\}$ are nested for each value of $a$, which increases in increments of four for a fixed color, it is clear that no two edges of the same color cross. 
As $a$ takes on all values from $1$ to $4k$, these pairs of nested edges cover $8k$ distinct edges of the graph. In the case where $a = 4k+1 = n$, the edge pair $\{a{-}(a{+}3), (a{+}1){-}(a{+}2)\}$ is $\{n{-}3, 1{-}2\}$, reducing mod $n$; these two edges fit on the fifth page without~crossing. 

This process, illustrated with $k = 3$ for the graph $C(13, \{1,3\})$ in Fig. 4, may be viewed as taking the red coloring and rotating it three more times, switching colors for each rotation. Now add a fifth color for the last two edges.

One may also describe this family of colorings, 
using  two Hamiltonian paths, alternatingly colored Red/Purple and Blue/Green.  In Fig. 4, left, the Red/Purple path is $1,4,3,2,5,8,7,6,
9,12,11,10,13$. 
As the reader will observe, there is a simple algorithm for the numbering: up 3, down 1, down 1, u3, u3, d1, d1, u3, u3, d1, d1, u3, and the resulting values mod $4$ are periodic: $1, 0, 3, 2, \ldots $.
The Blue/Green path is similar. The Hamiltonian paths each have $n-1$ edges; the two missing edges in black are the sparse page.
\begin{figure}[ht!]
\centering
\includegraphics[width=60mm]{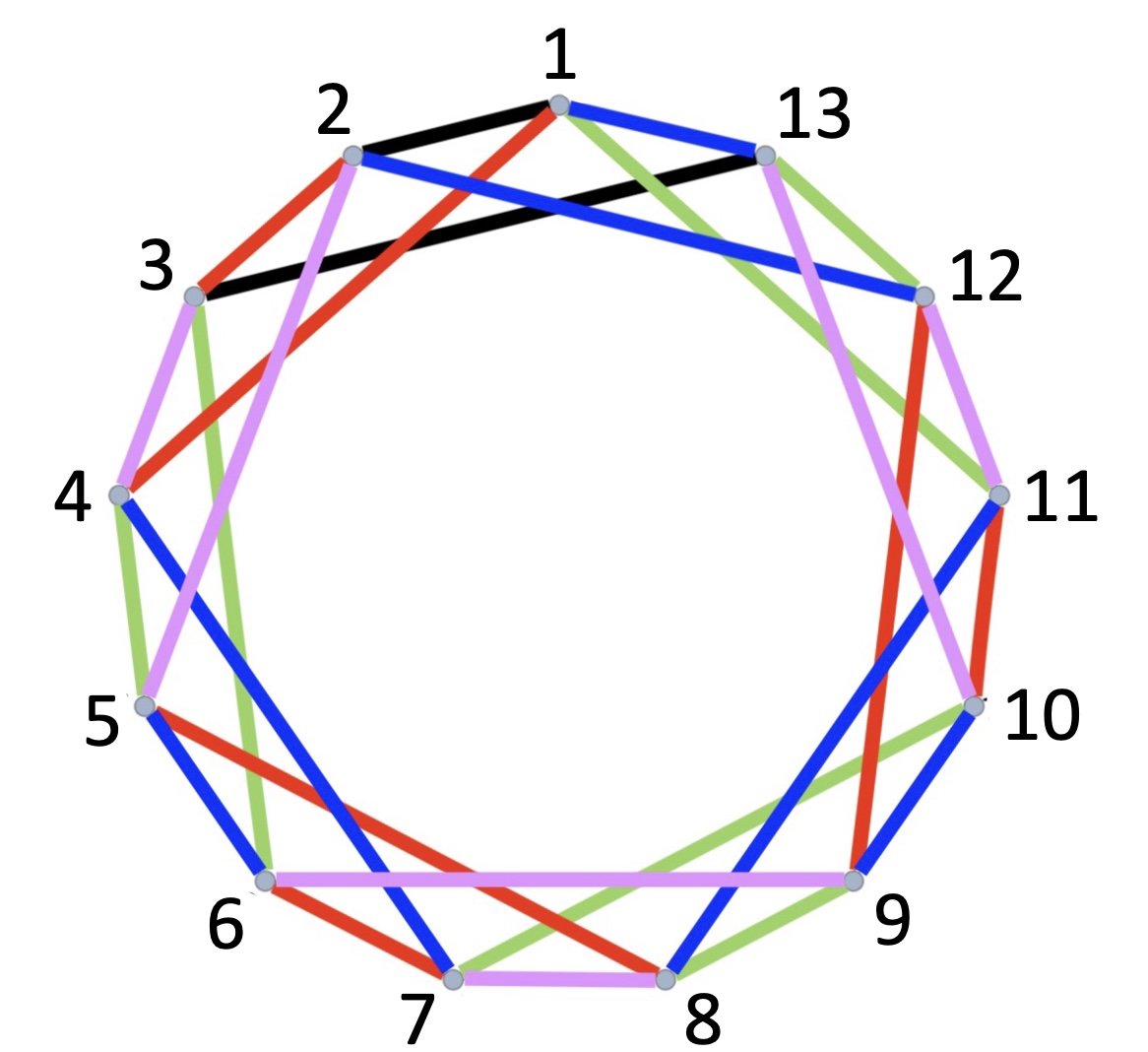}
\includegraphics[width=55mm]{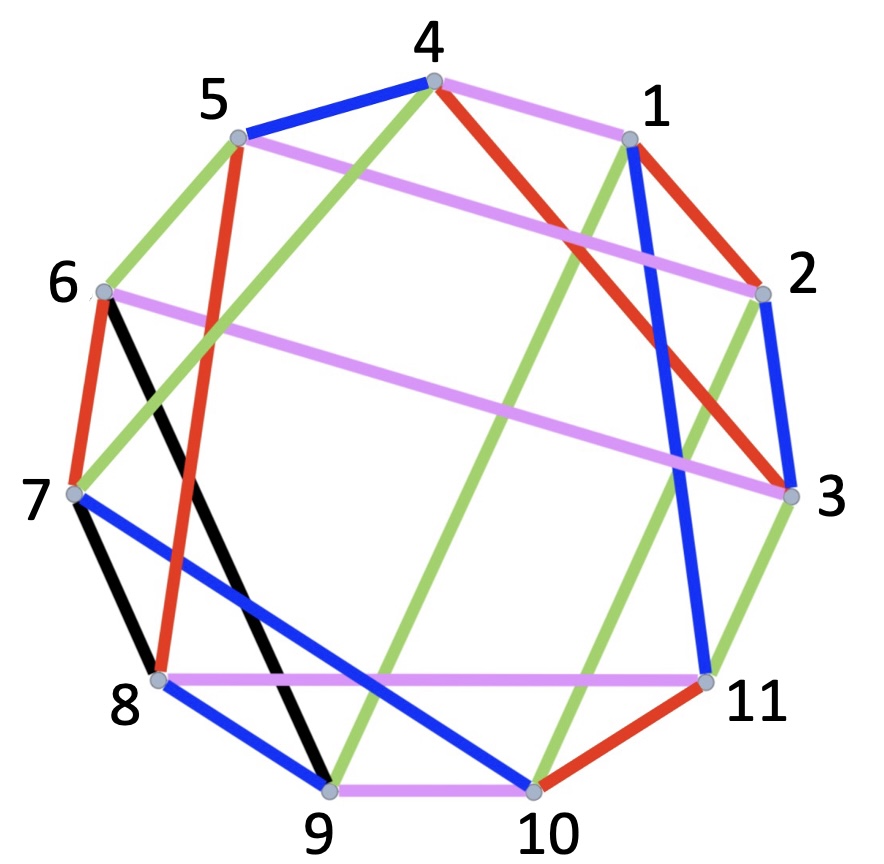}
\caption{Nearly dispersable embeddings of $C(13,\{1,3\})$ and $C(11,\{1,3\})$ \label{Figure 4}}
\end{figure}

{\bf Case 2}:  $n=4k+3$. We 
show that $mbt(C(n,\{1,3\}),\tau_n) = 5$, where $\tau_n$ is as in Theorem \ref{th:n-1-3-even} above.  Let $c$ be the same edge-coloring for the $10$ twist edges given in the  proof of Theorem \ref{th:n-1-3-even}.
If $k=1$ (so $n = 7$), color the non-twist edges $6{-}7$ red, $4{-}5$ blue, and  $4{-}7$ and $ 5{-}6$ black to get a 5-page matching book embedding.
If $k \geq 2$, then we have a general coloring scheme
\begin{itemize}
\item Green: $a{-}(a{+}3)$, $(a{+}1){-}(a{+}2)$, $a = 4t$, $t \in [k-1]$,
\item Red: $a{-}(a{+}3)$, $(a{+}1){-}(a{+}2)$, $a =  1+4t$, $t \in [k-1]$, and $n{-}(n{-}1)$,
\item Purple: $n{-}(n{-}3)$, $(n{-}1){-}(n{-}2)$ (Sparse; has 5 edges including twist),
\item Black: $a{-}(a{+}3)$, $(a{+}1){-}(a{+}2)$, $a =  2 + 4t$, $t \in [k-1]$,
\item Blue: $a{-}(a{+}3)$, $(a{+}1){-}(a{+}2)$, $a = 3+4t$, $t \in [k-1]$, and $4{-}5$.
\end{itemize}
Non-twist edges of the same color are either
nested pairs 
or 
or join consecutive vertices in the layout, so they cannot cross. The red and blue pages each contain $2(k-1) + 1 = 2k-1$ non-twist edges, the green and black pages each have $2(k-1) = 2k-2$ non-twist edges, and the purple page contains $2$ non-twist edges. The $5$-page matching book embedding covers all $10+2(2k-1) + 2(2k-2) + 2 = 8k + 6 = 2(4k + 3) = 2n$ distinct edges of the graph. Fig. 4 (right)  illustrates this coloring scheme for $C(11, \{1,3\})$ with $k = 2$.
\end{proof}

Note that in the above theorem, when $n > 7$ and $n = 4k + 3$, the two purple common twist edges $1{-}4$ and $2{-}5$ could be assigned the color black. This would reduce the number of purple edges from five to three. As $n \equiv 3$ (mod 4) increases, the sparse purple page stays at three edges. Hence, sparseness is 3, one more than the minimum.

\section{$C(n,\{1,2\})$ is nearly dispersable}

Now we consider the circulant graphs  $C(n, \{1,2\})$ for all $n \geq 4$. In the case $n = 4$ the corresponding graph is $K_4$, which has $mbt(G) = \Delta(G)+1 = 4$. 

\begin{theorem}
Let $n \geq 5$. Then $C(n,\{1,2\})$ is nearly dispersable.
\label{th:cn12}
\end{theorem}

\begin{proof}
We color this circulant in three cases using paths and cycles.

For $n \geq 5$ odd, draw the circulant using the odd-up, even-down cyclic order
\[\omega_n:=
(1, 3, \ldots, n, n{-}1, n{-}3, \ldots, 2).
\]
Color the edges of the ``zig-zag'' (Hamiltonian) path $1, 2, 3, \ldots, n$
using the colors red, blue alternatingly for the edges $(k{-}1){-}k$, where $k = 2, 3, \ldots, n$, along the path.   Color alternatingly with purple and green, the ``cross path'' 
$$3, 5, \ldots, n, 1, n{-}1, n{-}3, \ldots, 2.$$  This leaves two edges $1{-}3$ and $2{-}n$ each of which can be colored black. Thus, $mbt(G, \omega_n) = 5$ and the sparseness is minimum. See Fig. 5, left.
\begin{figure}[ht!]
\centering
\includegraphics
[width=45mm]{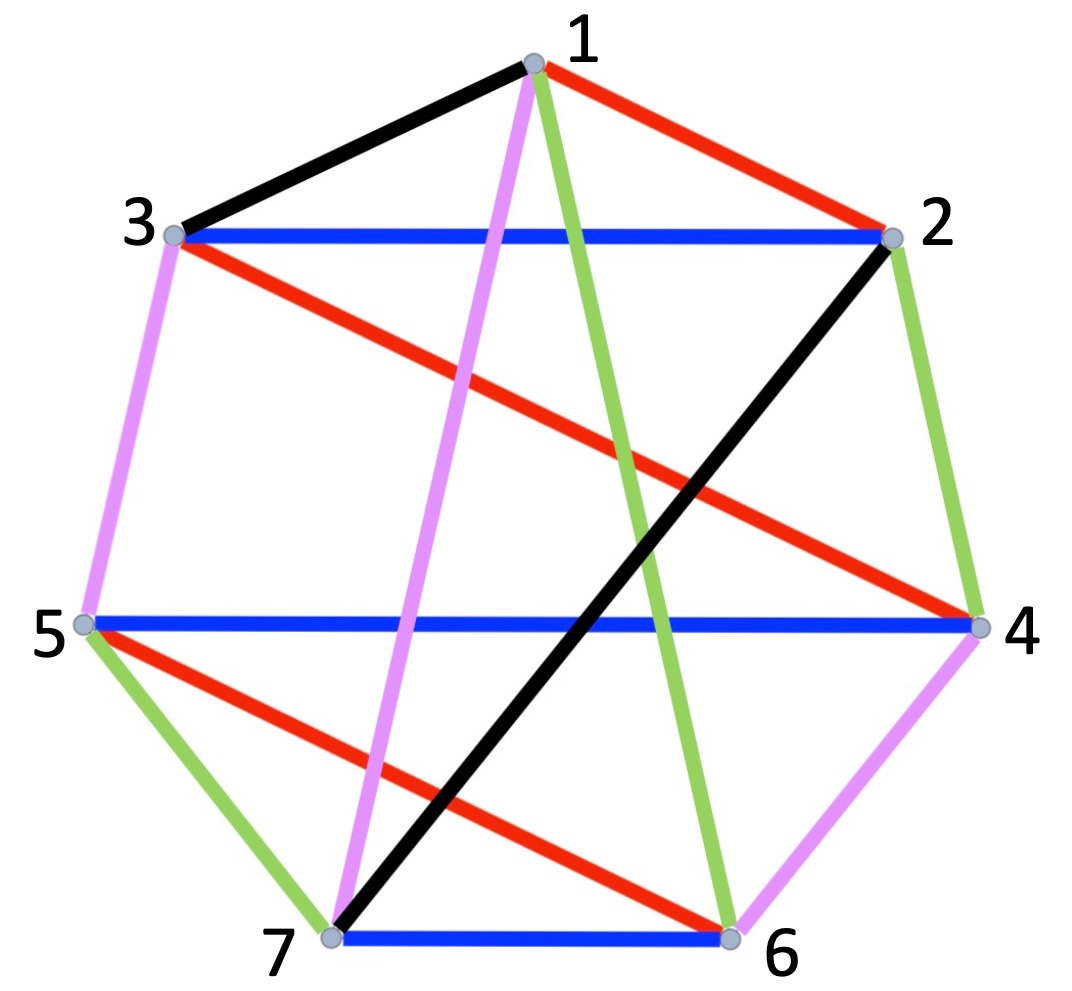}
\includegraphics
[width=45mm]{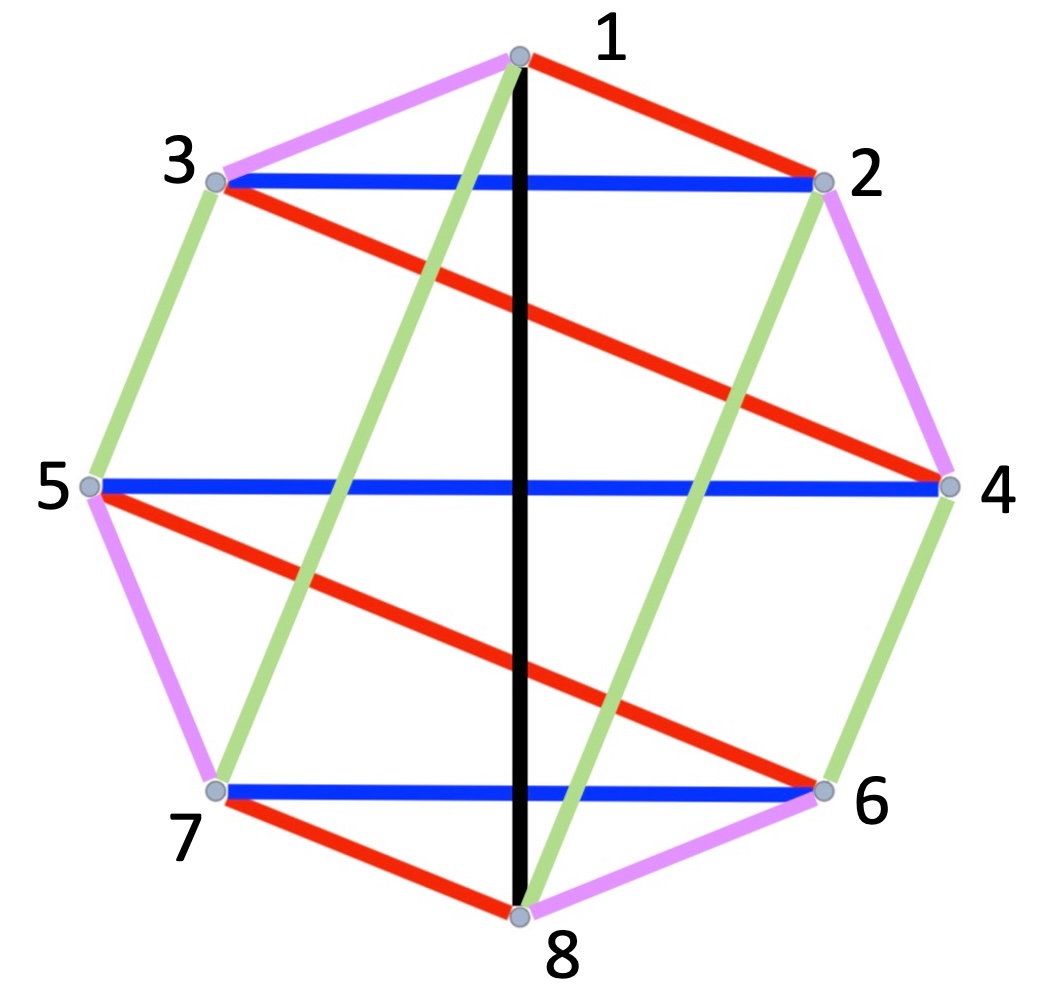}
\includegraphics
[width=45mm]{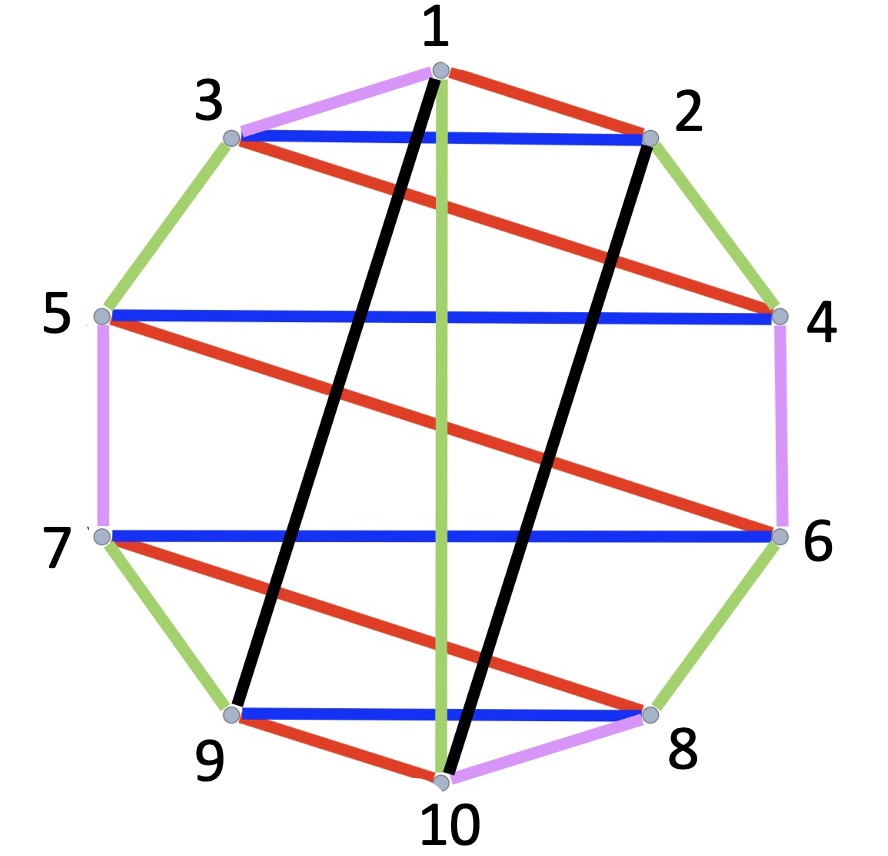}
\caption{$C(7,\{1,2\})$, $C(8,\{1,2\})$ $C(10,\{1,2\})$ are  nearly dispersable.\label{Figure 5}}
\end{figure}

For $n \geq 6$ even,  again use the odd-up, even-down cyclic order
\[
\omega_n := (1, 3, \ldots, n{-}1, n, n{-}2, \ldots, 2).
\]

As above, color the edges of the zig-zag path $1, 2, 3, \ldots, n$ using the colors red, blue alternatingly.    
If $n \equiv 0$ (mod 4), color the two disjoint non-self-crossing  cycles $ n{-}1, n{-}3, \ldots, 1, n{-}1$ and $2, 4, \ldots, n, 2$ (both of  length $n/2$) alternating purple and green. This accounts for $2n-1$ edges and the remaining edge $1{-}n$ is colored black. The sparseness is $1$ which is the minimum possible for an even order, nearly dispersable graph. See Fig. 5, center.

If $n \equiv 2$ (mod 4), edge-disjoint from the zig-zag path, there is another Hamiltonian path 
$n{-}1, n{-}3, \ldots, 
1, n, n{-}2, \ldots, 2$ and we color it alternatingly green and purple.  The two remaining edges $1{-}(n{-}1)$ and $2{-}n$ are parallel and colored black, see Fig. 5, right; they constitute the sparse page.
\end{proof}
\section{The case $C(n,\{2,3\})$}
\begin{figure}[ht!]
\centering
\includegraphics[width=45mm]{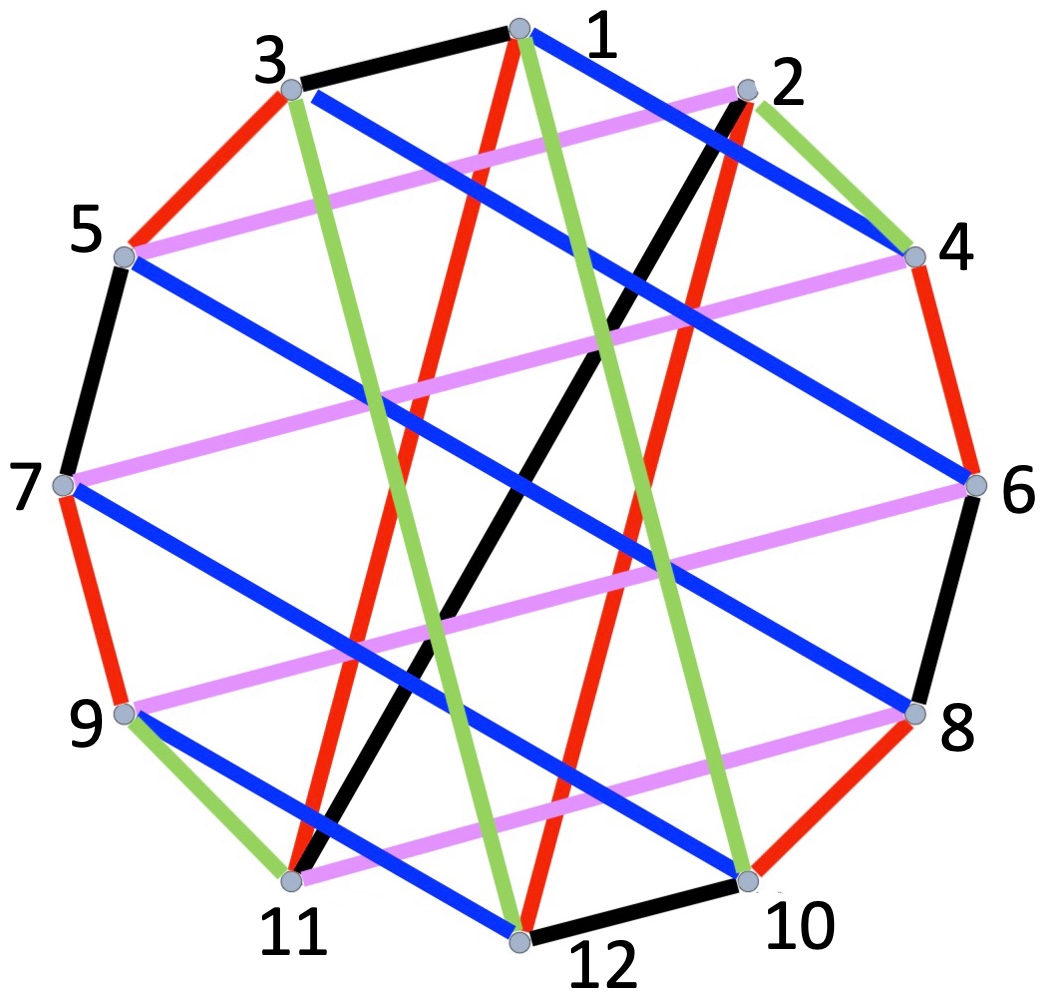}
\includegraphics[width=45mm]{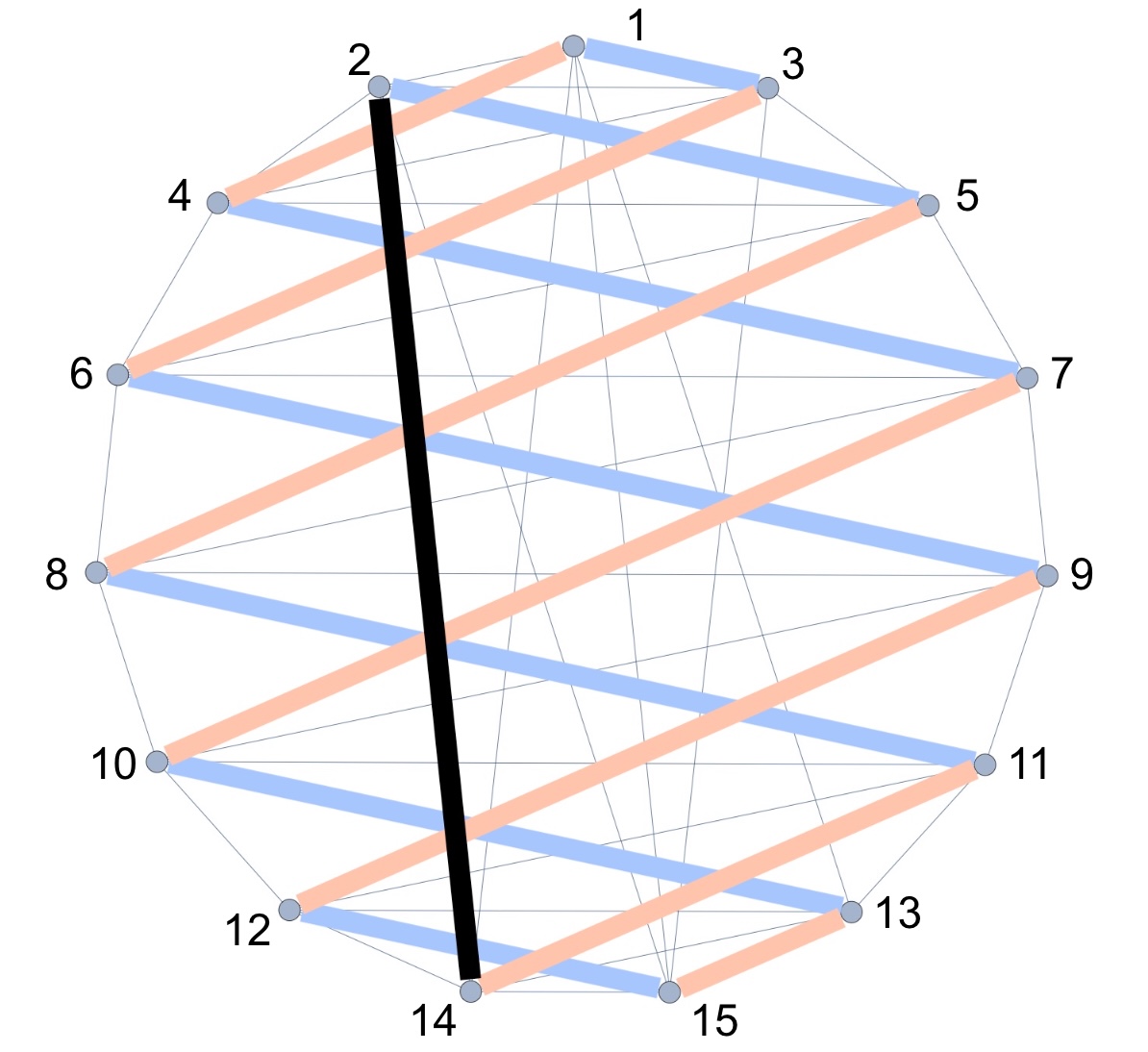}
\includegraphics[width=45mm]{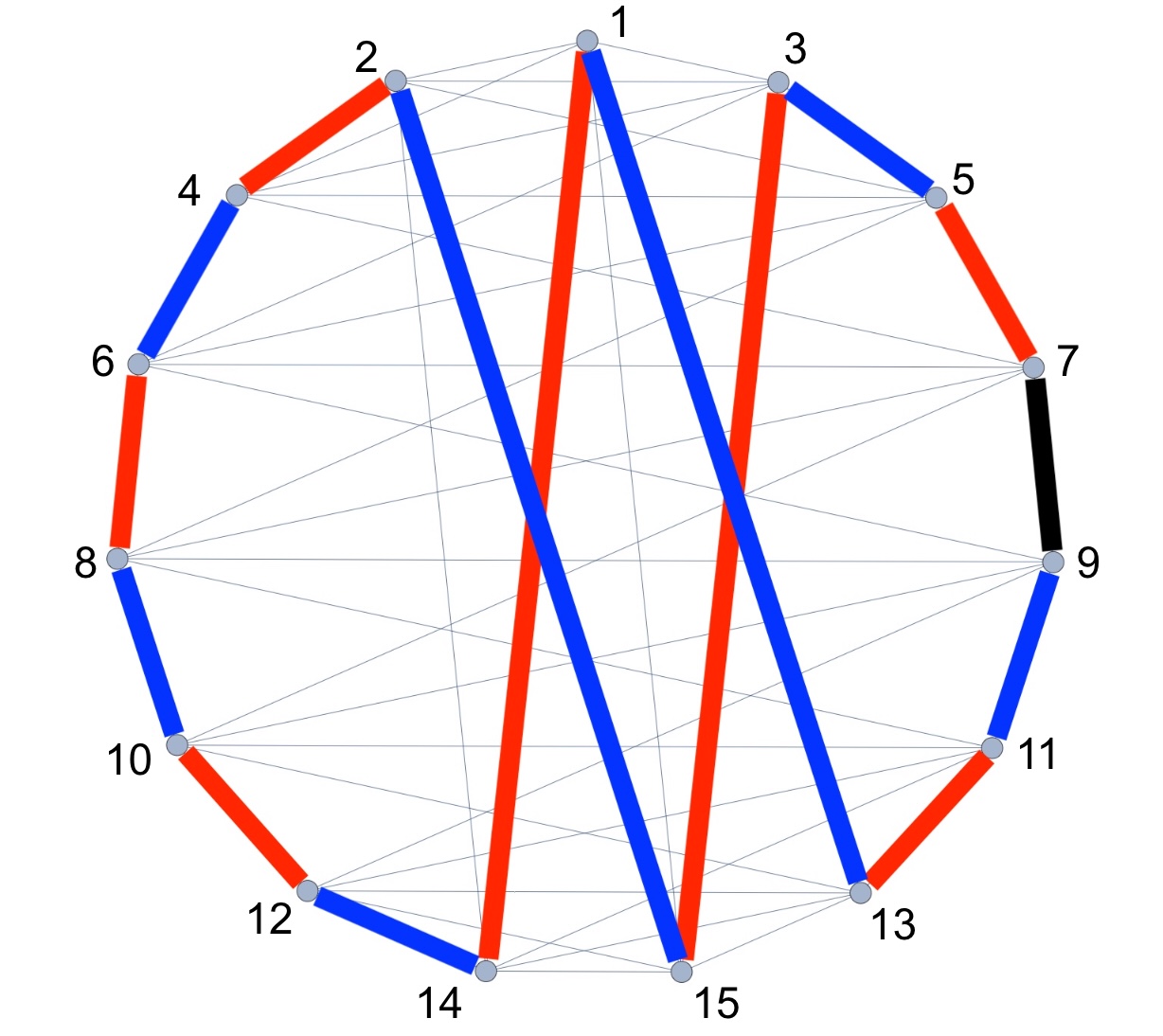}
\caption{$C(12,\{2,3\})$ on left, pages of $C(15,\{2,3\})$ in middle and right. \label{Fig:12-2-3}}
\end{figure}

Let $\omega_n$ be the odd-up, even-down vertex-order from the proof of Theorem \ref{th:cn12}. 
\begin{theorem}
If $n \geq 8$ is {\rm even}, then $(C(n,\{2,3\}),\omega_n, c_n)$ is nearly dispersable, where $c_n$ is the coloring given below.
\label{th:cn23-even}
\end{theorem}
\begin{proof}
Put $n=2k$; consider the Hamiltonian path in $C(n,\{2,3\})$ given by 
\[
(n{-}3,n{-}5,\ldots,1,n{-}1,2,n,n{-}2,\ldots,4)
\]
which is noncrossing w.r.t. the layout $\omega_n$; we color it alternatingly red and black.  This colors $2k-1$ edges. Color with purple the pairwise-disjoint edges $2{-}5$, $4{-}7$, $\ldots$, $(n{-}4){-}(n{-}1)$, which accounts for $k-2$ edges. Color with blue the pairwise-disjoint edges $1{-}4$, $3{-}6$, $\ldots$, $(n{-}3){-}n$, accounting for another $k-1$ edges.  In the sparse green page, there are four edges: $2{-}4$, $(n{-}3){-}(n{-}1)$, $3{-}n$, and $1{-}(n{-}2)$. Hence, all $4k$ edges are used. 
As same-color edges do not cross or share an endpoint,
$\omega_n$ is nearly dispersable. See Fig. \ref{Fig:12-2-3}, left.
\end{proof}

That $C(6,\{2,3\})$ is nearly dispersable is left to the reader.  

\begin{theorem}
For $n \geq 7$ {\rm odd}, $C(n,\{2,3\})$ is nearly dispersable w.r.t. $\omega_n$.
\label{th:cn23-odd}
\end{theorem}
\begin{proof}
Let $n = 2r+1 \geq 7$. We decompose the edge set into two edge-disjoint Hamiltonian cycles. These odd-length cycles are 2-colored except for one black edge, placed on the sparse  page. 
The first of these two cycles, see Fig. \ref{Fig:12-2-3}-middle, is colored orange, aqua, and black according to the following scheme:
\begin{itemize}
\item Orange: $a{-}(a{+}3)$, $a = 1 +  2t$, $t \in [r-1]-1$ and $n{-}(n{-}2)$,
\item Aqua: $a{-}(a{+}3)$, $a =  2t$, $t \in [r-1]$ and $1{-}3$,
\item Black: $2{-}(n{-}1)$.
\end{itemize}
The second of these cycles, illustrated in Fig. \ref{Fig:12-2-3}-right, is colored red, blue, and black according to the following scheme:
\begin{itemize}
\item Red: $1{-}(n{-}1)$, $3{-}n$, and $r{-}2$ additional edges of the form $a{-}(a{+}2)$ on the outer cycle (alternating with the blue edges), 
\item Blue: $2{-}n$, $1{-}(n{-}2)$, and $r{-}2$ additional edges of the form $a{-}(a{+}2)$ on the outer cycle (alternating with the red edges),
\item Black: $r{-}(r{+}2)$.
\end{itemize}
The orange and aqua pages each consist of a total of $r$ non-crossing parallel edges. The red and blue pages also each contain $r$ non-crossing edges with two parallel edges on each page through the center and the remaining $r{-}2$ edges on the outer boundary. The remaining two black edges, which is the minimum possible number, clearly do not intersect on the sparse page since the $r{-}(r{+}2)$ edge lies on the outer cycle and does not share an endpoint with $2{-}(n{-}1)$. Hence, all $4r + 2 = 2n$ edges are accounted for.\end{proof}

\begin{figure}[ht!]
\centering
\includegraphics[width=60mm]{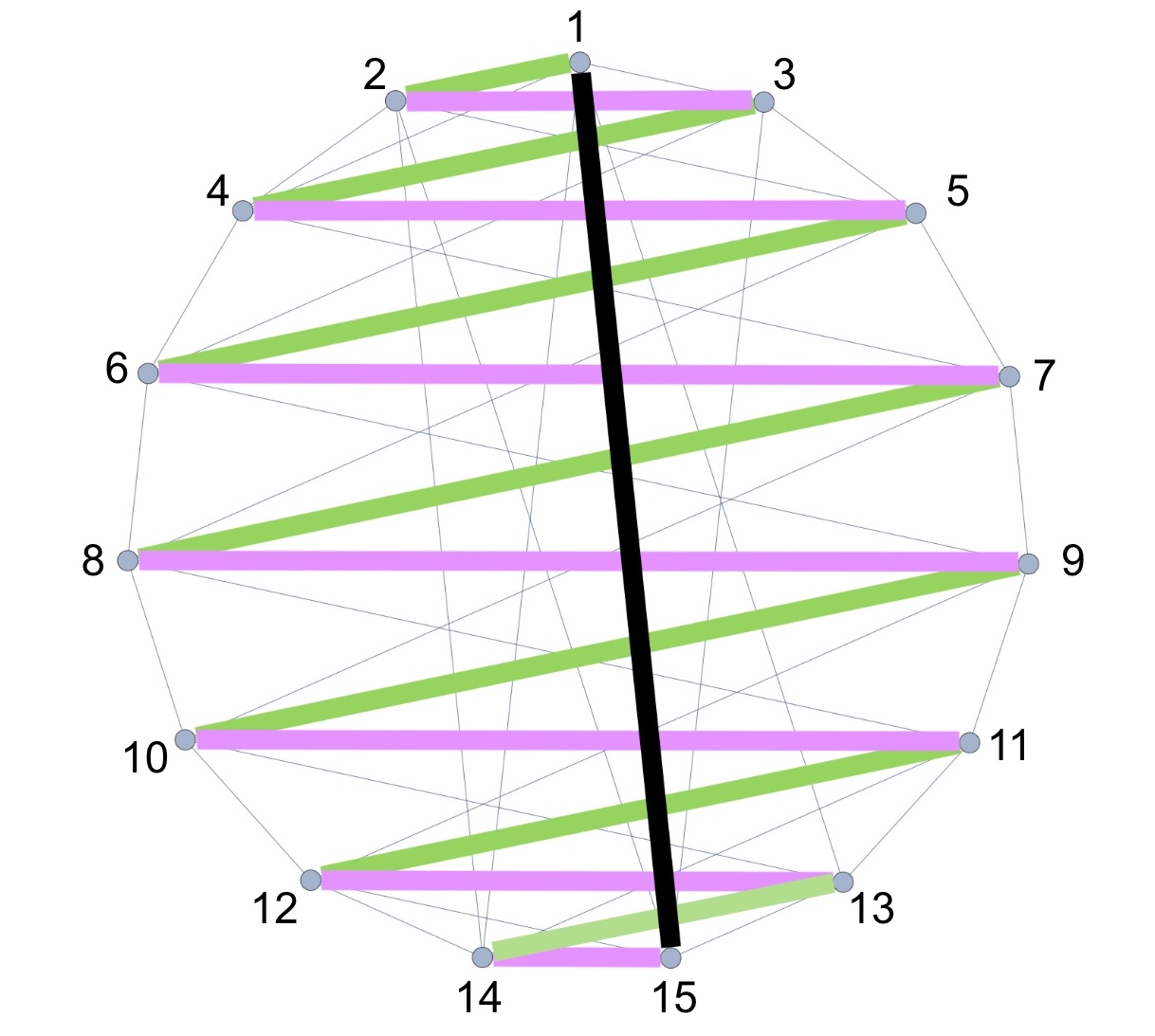}
\includegraphics[width=60mm]{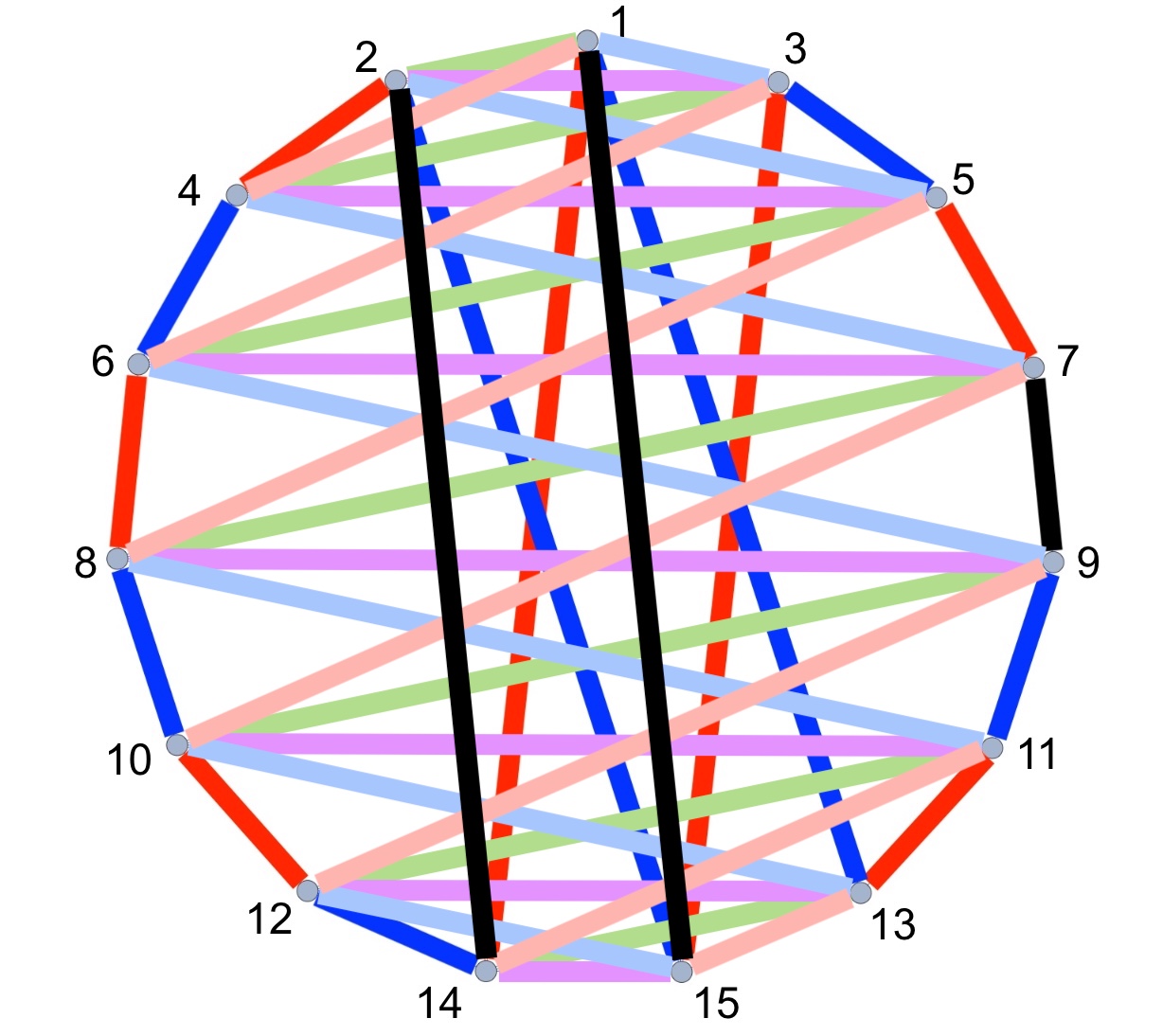}
\caption{ Length-1 edges (left) and 7-coloring of $C(15,\{1,2,3\})$ (right).\label{fig-superimp}}
\end{figure}

\begin{theorem}
For $n \geq 7$ {\rm odd}, $C(n,\{1,2,3\})$ is nearly dispersable w.r.t. $\omega_n$.
\end{theorem}
\label{th:cn123-odd}
\begin{proof}

Let $n = 2r+1 \geq 7$. Use the identical layout and coloring scheme as in Theorem 5. This will cover all distance 2 and distance 3 edges. Now observe that a purple-green non-crossing Hamiltonian path $1, 2, 3, \ldots, n$ can be added to cover all of the distance-1 edges, with the exception of edge $1{-}n$. This last edge can be placed on the black (sparse) page and does not intersect either $r{-}(r{+}2)$ or $2{-}(n{-}1)$. We note that this new cycle contributes $r$ purple edges, $r$ green edges, and $1$ black edge as shown in Fig. \ref{fig-superimp} left. Combining this with $C(n,\{2,3\})$, we have accounted for all $6r+3 = 3n$ edges of $C(n,\{1, 2,3\})$ and have achieved an optimal sparseness of $3 = \Delta/2$. See Fig. \ref{fig-superimp} right for the combined nearly-dispersable coloring of $C(n,\{1,2,3\})$ for $n$ odd.
\end{proof}

It is natural to ask if the above scheme for adding distance-1 edges allows {\it even values} of $n$.  This almost works, with the exception of the black edge $1{-}n$, which intersects edges on the sparse page in the above layout for even values of $n$, so for even $n$, $mbt(C(n,\{1, 2,3\})-e) = 7$, where $e=1{-}n$.

\section{Larger degree and jump-lengths}

We now give minimum layouts of $C(n,\{1, 2,3\})$ for some even values of $n$ and for a variety of circulants with $\Delta > 3$, using polymerization and periodicity.

\begin{theorem}
With  $\nu_n = (1,2,\ldots,n)$ and $r=2k+1 \geq 5$, suppose
$r | n$. Then $$mbt(C(n,\{1,2,\ldots,k\}),\nu_n) = 2k+1.$$
\label{th:n[k]}
\end{theorem}
\vspace{-1cm}

\begin{figure}[ht!]
\centering
\includegraphics[width=110mm]{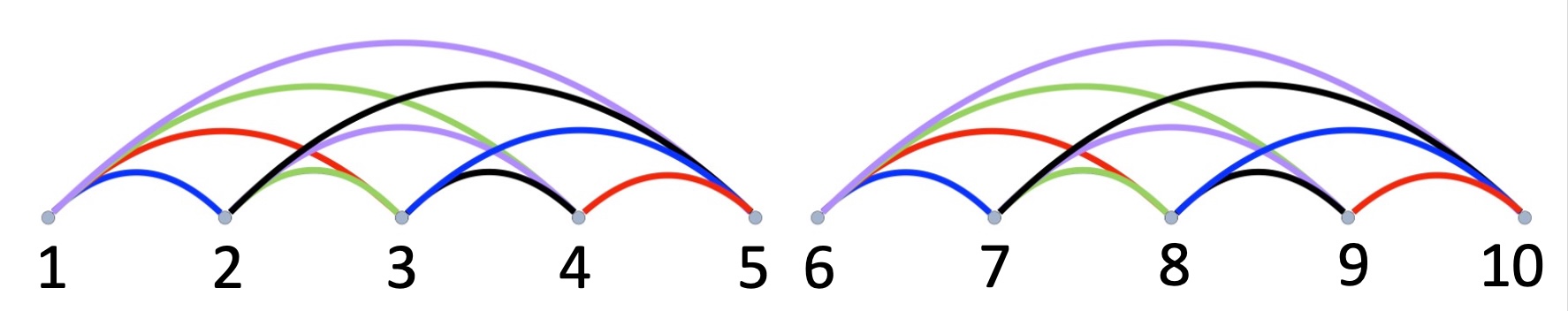}
\caption{Two copies of $K_5$ . }
\label{Figure 7}
\end{figure}

\begin{figure}[ht!]
\centering
\includegraphics[width=100mm]{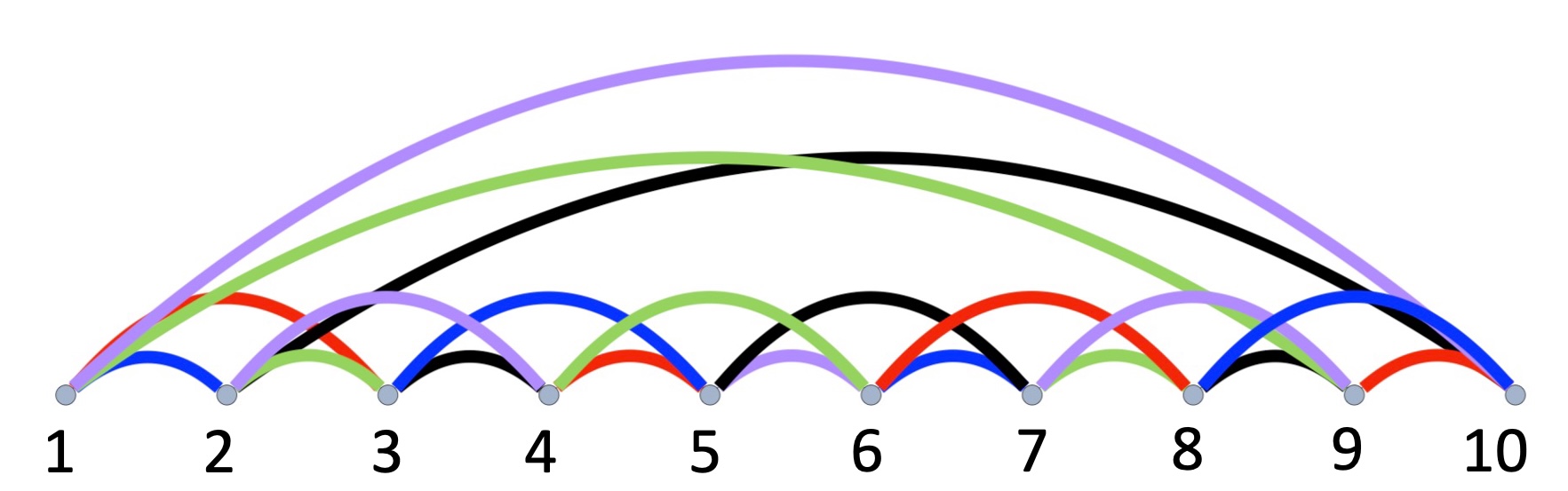}
\caption{$C(10,\{1,2\})$ obtained by polymerization.
\label{Figure 8}}
\end{figure}

\begin{proof}
Use Lemma \ref{lm:polymer} on the embedding  $(K_r, \nu_r, c)$ in \cite{so-thesis}; see Fig \ref{Figure 7} and \ref{Figure 8}.
\end{proof}


For $r=7$,
and for all $m\geq 1$, one has the consequence (new for $m$ even):
\begin{equation}
mbt(C(7m, \{1,2,3\})) = 7.
\label{eq:7m}
\end{equation}
Similarly, $mbt(C(9m, \{1,2,3,4\})) = 9$, $mbt(C(11m, \{1,2,3,4,5\})) = 11$, etc.

Note that, for $k$ a positive integer, $C(n,[k]) \cong K_n$ for $n \in \{2k, 2k+1\}$, but the isomorphism determines {\it which} edge-set $E_j$ of $K_n$ corresponds to the length-$j$ jumps in the circulant for $j = 1, \ldots, k$.  With $n = 2k$ or $2k+1$, the length-$k$ jumps are longest in the circulant and they induce a 1-factor or a spanning cycle according to whether $n$ is even or odd.

Corresponding to the $n=2k$ case,
for $k \geq 3$, the {\bf cocktail party graph} $O_k := K_{2k} -k K_2$ is the complement of a 1-factor. It is also the 1-skeleton of the $n$-dimensional {\bf octahedron} and is regular with $\Delta O_k = 2k-2$.  As it contains triangles, the octahedron is at best nearly dispersable.  We show that it does have a nearly dispersable embedding, with $2k-1$ pages, but the embedding does not use the standard vertex ordering and so we don't have a direct way to polymerize it.
However, the natural vertex order gives a matching book embedding with one additional page and this {\it can} be polymerized.  
For the $n=2k+1$ case, the same things can be done with $\overline{C_n}$, the complement of $C_n$.

The {\bf folded order} $\phi$ from \cite{pck-circ} will be needed for $n \in \{2k, 2k+1\}$.
$$\phi_n := (1, 2, \ldots, k, n, n-1, \ldots, k+1).$$

\begin{figure}[ht!]
\centering
\includegraphics[width=60mm]{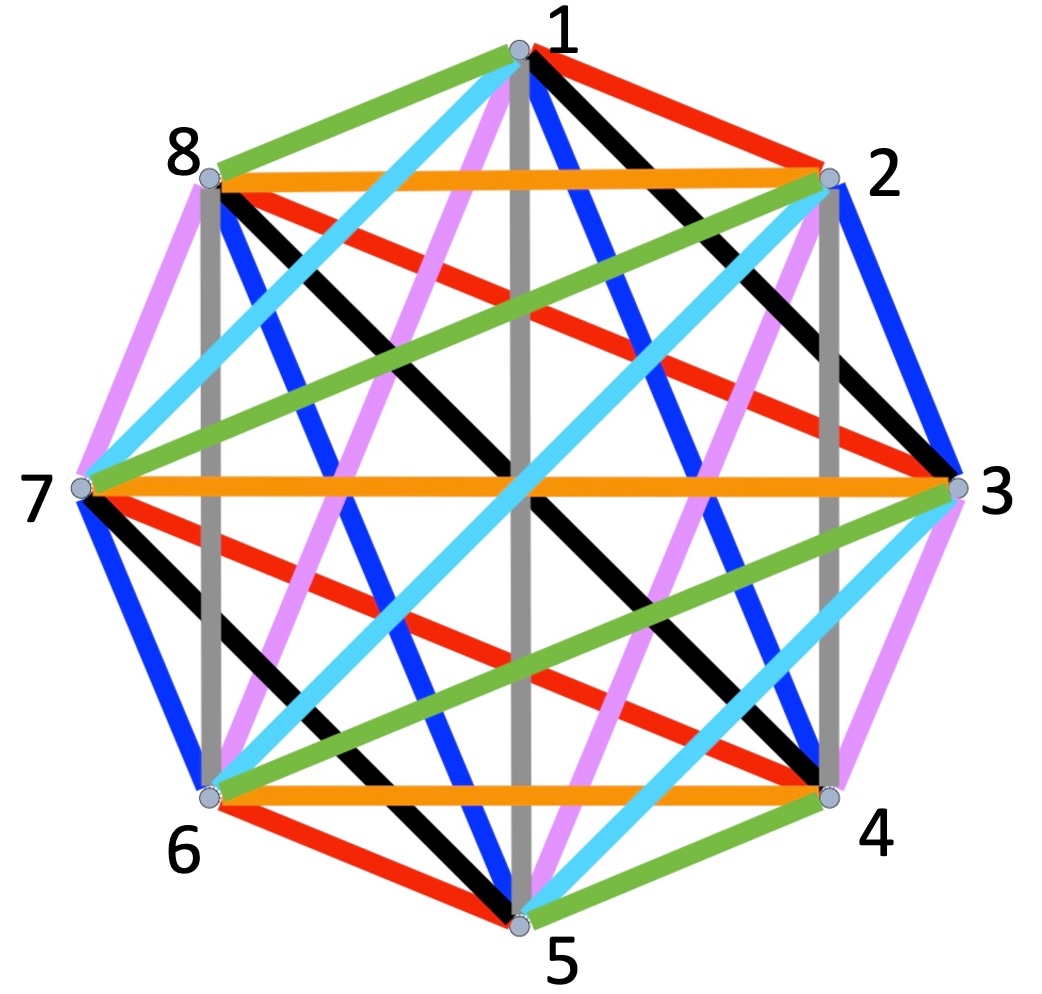}
\caption{$K_8$. \label{Fig:8-1-2-3}}
\end{figure}

\begin{figure}[ht!]
\centering
\includegraphics[width=55mm]{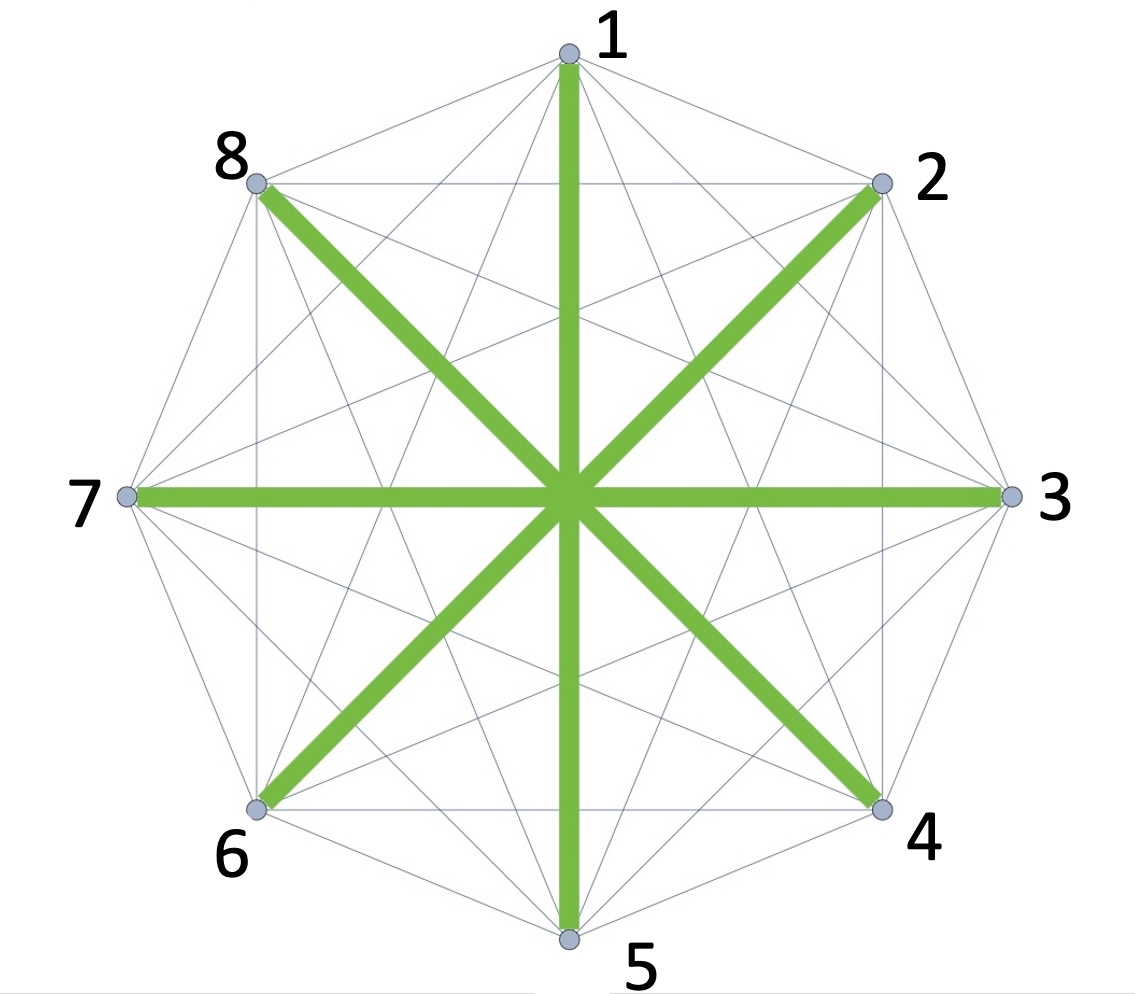}
\includegraphics[width=55mm]{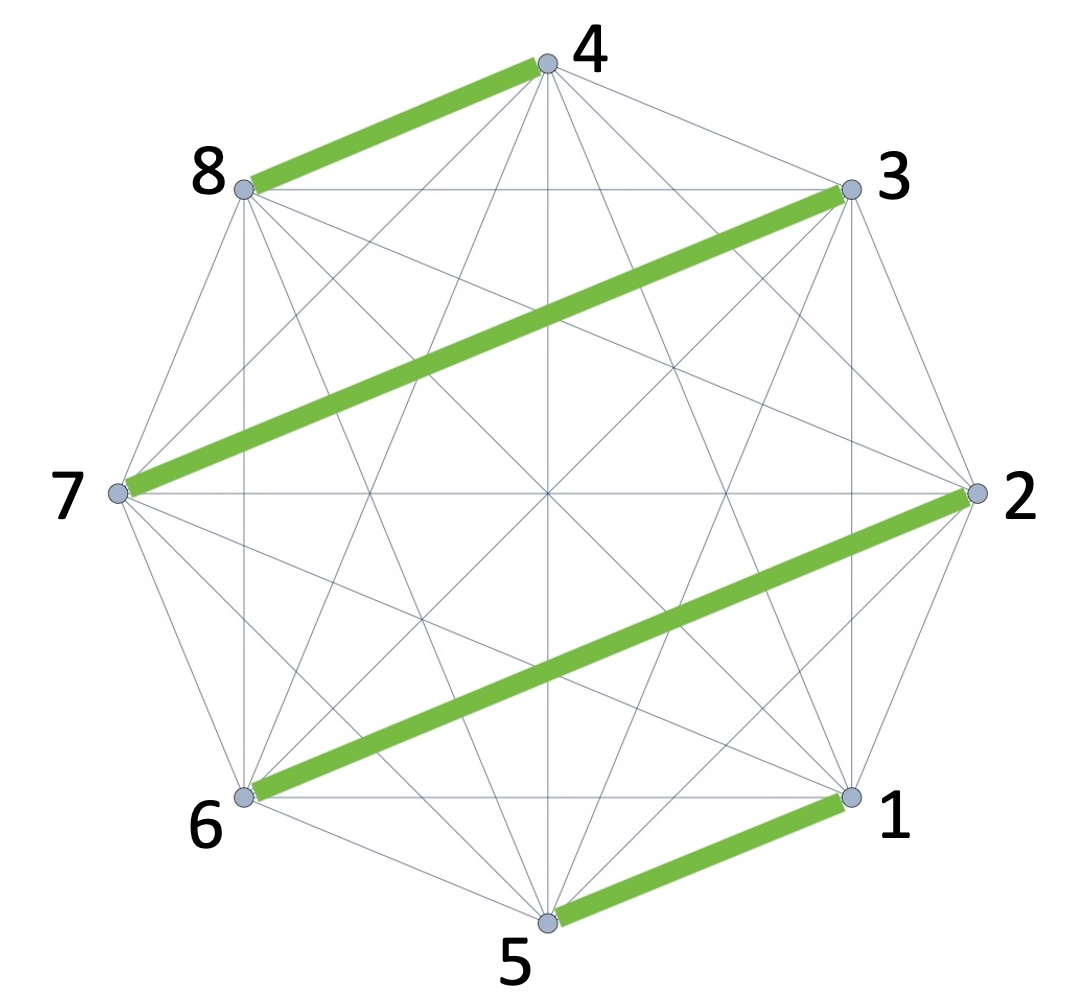}
\caption{A perfect matching in $K_8$. \label{Fig:8-m}}
\end{figure}

\begin{figure}[ht!]
\centering
\includegraphics[width=70mm]{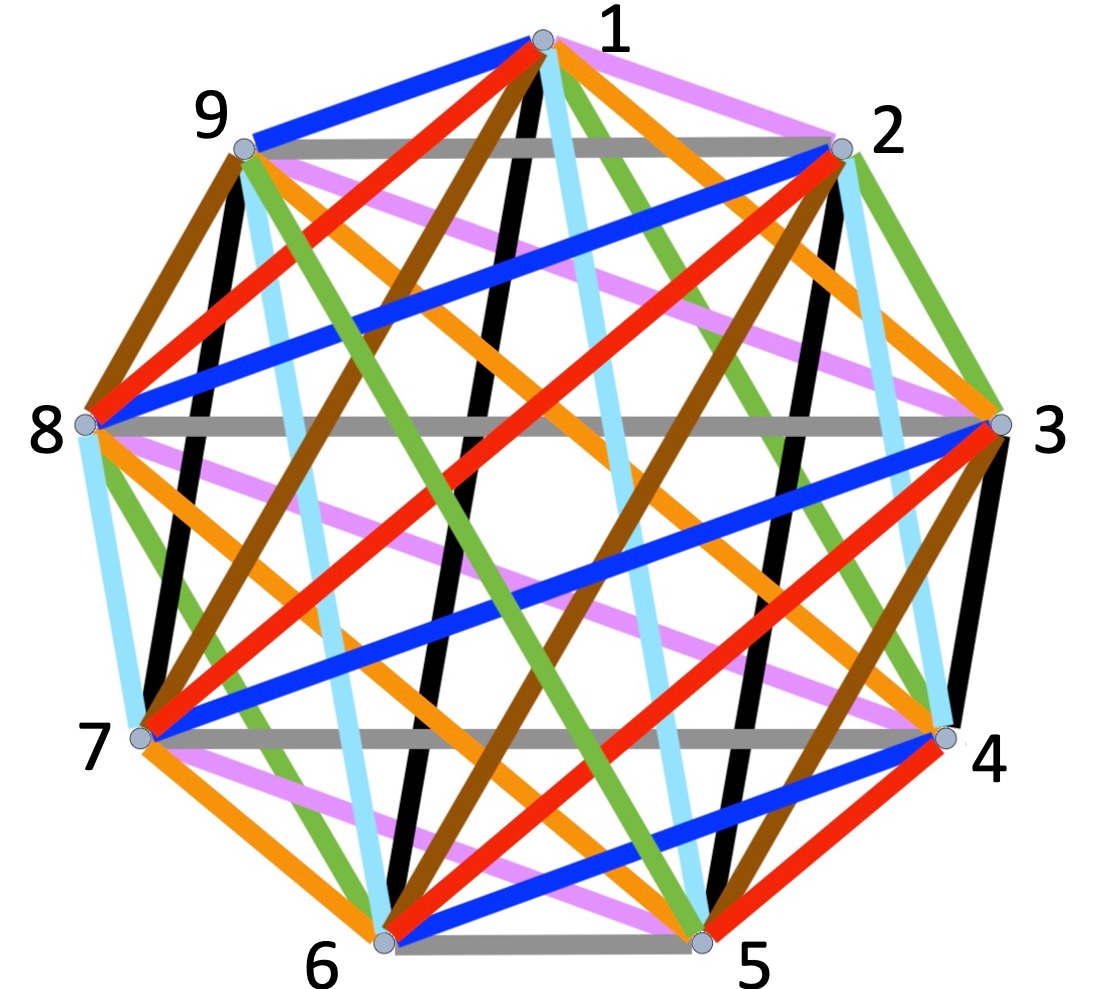}
\caption{$K_9$. \label{Fig:9-123}}
\end{figure}

\begin{figure}[ht!]
\centering
\includegraphics[width=60mm]{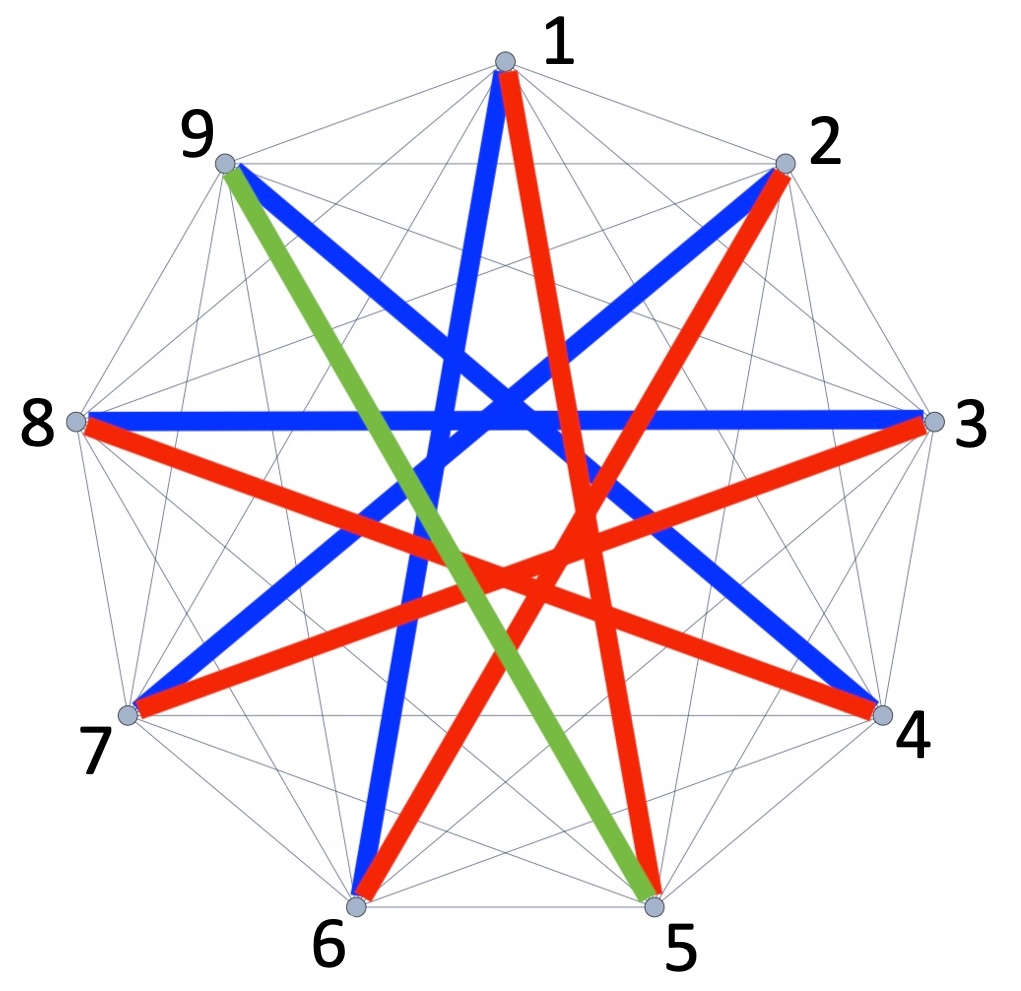}
\includegraphics[width=60mm]{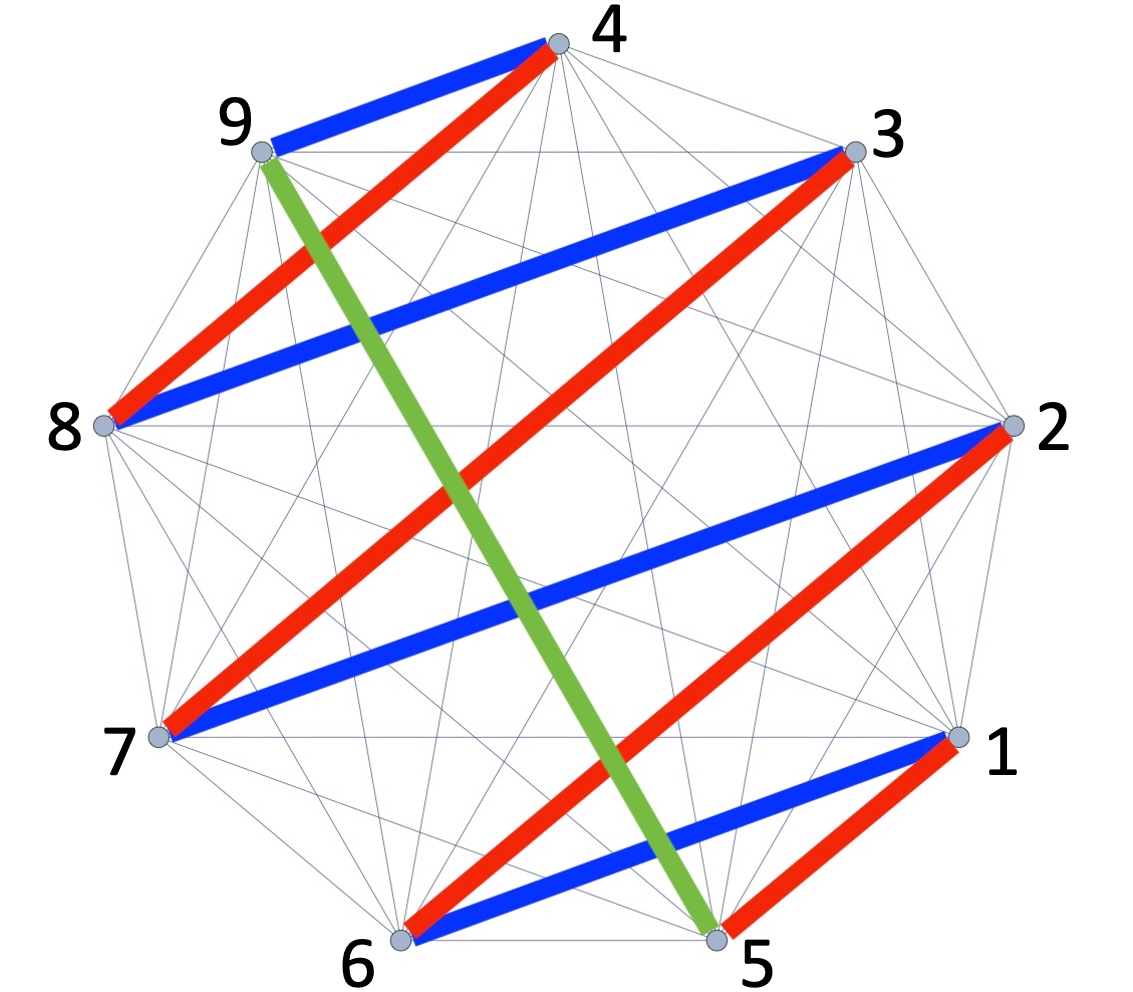}
\caption{A Hamiltonian cycle in $K_9$. \label{Fig:9-star}}
\end{figure}

\begin{theorem}
For $k \geq 3$, $O_k$ and $\overline{C}_{2k+1}$ are nearly dispersable with
\begin{equation}
mbt(O_k, \phi_{2k}) = 2k-1 \,=\, mbt(\overline{C}_{2k+1}, \phi_{2k+1}).
\label{eq:8m}
\end{equation}

\label{th:S=[k-1]}
\end{theorem}
\begin{proof}
For $n \geq 6$, even or odd, the complete graph $K_n$ is nearly dispersable using the natural vertex order with pages being the 1-factors produced by maximal families of parallel edges given in \cite[p 87]{so-thesis}; see Fig. \ref{Fig:8-1-2-3} and Fig. \ref{Fig:9-123}.

By our remark above about how the edges can correspond to various jump-lengths in the isomorphic circulant graph, we note that for $n=2k$ and the folded order $\phi_{2k}$,
the set of length-$k$ edges, $E_k$, are a parallel matching, one of the $2k$ pages in the nearly dispersable matching book embedding of $K_{2k}$, and the removal of $E_k$ leaves a $2k-1$-page layout of $O_k$.  See right side of Fig. \ref{Fig:8-m}.

For $n=2k+1$, the length-$k$ edges in the circulant constitute a Hamiltonian cycle $Z$.  The folded order takes this $2k+1$-cycle into 2 of the $2k+1$ pages of the standard matching book embedding of $K_{2k+1}$ given in \cite{so-thesis}, with one extra edge. Deleting $Z$ gives a $2k-1$-page layout of $\overline{C}_{2k+1}$ (Fig. \ref{Fig:9-star}, right).
\end{proof}
With one page above the minimum, we have {\it polymerizable} embeddings,
$$mbt(O_k, \nu_{2k}) = 2k = mbt(\overline{C}_{2k+1}, \nu_{2k+1}).$$  Under the natural order, the edges of $K_{2k}$ which are length-$k$ edges in $C(2k, [k])$ pass through the center of the circle. 
Hence, removing them only decreases by 1 the number of edges in each of the $2k$ pages.  See left side of Fig. \ref{Fig:8-m}. Similarly, with the layout $\nu_{2k+1}$, the length-$k$ edges form a mandala-like figure (Fig \ref{Fig:9-star}, left) which can be deleted from the matching book embedding.

Thus, if $n$ is a multiple of $2k$ or of $2k+1$ for $k\geq 3$, then 
by Lemma \ref{lm:polymer}, 
\begin{equation}
mbt(C(n,\{1,\ldots,k-1\}), \nu_n) \leq 2k.
\label{eq:both}
\end{equation}


\medskip

\begin{figure}[ht!]
\centering
\includegraphics[width=60mm]{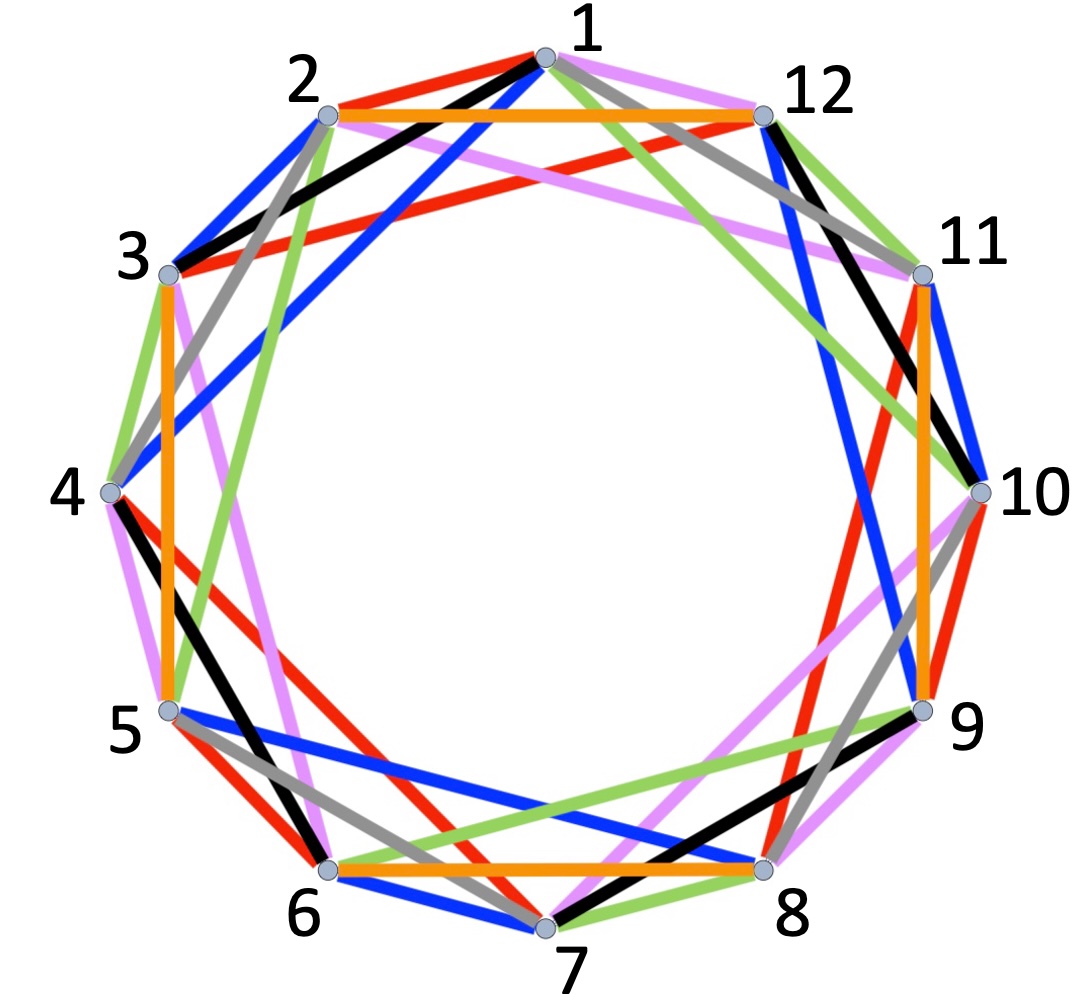}
\caption{$C(12,\{1,2,3\})$. \label{fig-12reg}}
\end{figure}

The next theorem allows multiples of 12.
\begin{theorem}
If $n=12m$, $m\geq 1$, then $C(n,\{1,2,3\})$ is nearly dispersable.
\label{th:12m}
\end{theorem}
\begin{proof}
An explicit coloring of the edges with respect to natural order $\nu_n$ using $7 =1+\Delta$ colors is given as follows: 
Periodically, 4-color the edges of length 1, and for every edge of  length 3, use the same color as the unique edge of length 1 with which it is nested, thus 4-coloring all edges of length 1 or 3.  
Three new colors (periodically) suffice for the remaining edges.  
See Fig. \ref{fig-12reg}.
\end{proof}

For the bipartite case, if $k \geq 1$, by Lemma \ref{lm:polymer}, 
analogous to (\ref{eq:both}), we have:

\begin{equation}
\mbox{If} \; 4k | n , \; \mbox{then} \;mbt(C(n,\{1, 3, 5, \ldots,2k-1\}), \nu_n) = 2k;
\label{eq:both1}
\end{equation}
\begin{equation}
\mbox{If} \; (4k+2) | n, \; \mbox{then} \;mbt(C(n,\{1, 3, 5, \ldots,2k-1\}), \nu_n) \leq 2k + 1.
\label{eq:both2}
\end{equation}

\section{Discussion}

Recently, Yu, Shao and Li \cite{yu-shao-li} have shown dispersability (or near dispersability) for circulants of degree 3 and degree 4  with  {\it all}  $\,$jump-lengths according to whether or not they are bipartite.  This extends our results for these degrees.
Our methods supply different solutions to the problem of finding optimal page-number matching book embeddings for such circulants.  

We have also considered higher degree circulants and have analyzed some of the structural features of matching book embeddings of regular graphs.

For $C(n,\{1,2,3\})$, we show near dispersability when $n$ is odd and 
when $n$ is even and divisible by $7$ or $12$. We also show near dispersability for the circulants resulting by deleting a maximum matching from an even-order $K_n$ and a spanning cycle from an odd-order $K_n$.  The analogous results holds in the bipartite case.

Results are constructive, not just existential, 
and so will remain useful even if the full conjecture on
vertex-transitive graphs is proved (or disproved).

Some condition on the graph is needed as a regular graph can have an {\it arbitrarily large} value for the ratio $mbt/\Delta$ according to Alam et al. \cite{abdgkp2021}, which uses a counting argument of McKay \cite{mckay} to prove that, for any fixed $\Delta \geq 3$, there exist $\Delta$-regular bipartite graphs $G_n$ with $mbt(G_n) \to \infty$ as $n \to \infty$. 

But vertex-transitive graphs are rather special. 
Du, Kutnar \& Maru\v{s}i\v{c}~\cite{du2021} showed that the Lovasz conjecture ({\it Every vertex transitive graph contains a Hamiltonian cycle, with five exceptional cases}) is correct when the order is a product of two primes and the graph satisfies additional conditions involving the action of a group.
Also, Diestel \cite[p 52; Ex 12]{diestel} notes that every connected, even-order, vertex-transitive graph has a 1-factor.

In many cases, our proof of near dispersability for a family of matching book embeddings uses a sparse page with the minimum number of edges. Indeed, the lower bound of $2$ is achieved in the proof of Theorem \ref{th:cn13-odd} for $C(n,\{1,3\})$ when $n \equiv 1$ (mod $4$) while we get sparseness $s \leq 3$ for  $n \equiv 3$ (mod $4$).
The proof of Theorem \ref{th:cn12} shows $C(n,\{1,2\})=2$ for $n$ odd, while for $n \equiv 0$ (mod $4$), sparseness $=1$ is achieved but for $n \equiv 2$ (mod $4$), we only get sparseness $\leq 2$. 
The proof of Theorem \ref{th:cn23-even} shows $s(C(n,\{2,3\})) \leq 4$ if $n$ is even.  Sparseness has minimum value (2 and 3, resp.) for Theorems \ref{th:cn23-odd} and \ref{th:cn123-odd} when $n$ odd and for degree 6.
Our matching book embeddings for 
Theorems \ref{th:n[k]}, \ref{th:S=[k-1]} and \ref{th:12m}, in contrast, are quite symmetric and so the opposite of sparse embeddings. 

For a nearly dispersable embedding,
sparseness allows the deletion of a small number of edges to eliminate an entire page, while symmetry might be preferred for an ``online'' problem where the graph being embedded is evolving.

Matching book thickness for regular graphs has a clear lower bound, so a coloring which achieves the minimum is detectable. We think that
finding the matching book thickness of various graphs could be a good target for genetic algorithms, neural networks, or artificial intelligence. 
See, for example, \cite{online}.
Application of machine learning techniques to matching book thickness might improve computational  theory and practice; see, e.g., \cite{hutson}. 
Indeed, one has an endless supply of vertex-transitive graphs on which to test procedures.

\subsection*{Acknowedgement}
We appreciate the referees for reading and re-reading our paper and for their helpful and constructive comments.

\end{document}